\documentclass[12pt]{article}
\usepackage{enumerate}
\usepackage{amsfonts}
\usepackage{amsmath}
\usepackage{amssymb}
\usepackage{amsthm}

\def\D{\mathcal{D}}
\def\V{\mathcal{V}}

\def\H{\mathcal{H}}
\def\K{\mathcal{K}}

\def\S{\mathfrak{S}}
\def\F{\mathfrak{F}}
\def\C{\mathfrak{C}}
\def\T{\mathfrak{T}}
\def\B{\mathfrak{B}}

\def\Y{\mathfrak{Y}}

\newcommand{\id}{\mathrm{Id}}
\newcommand{\Tr}{\mathrm{Tr}}

\newcommand{\shs}{\hspace{1pt}}
\newcommand{\sn}{\|\hspace{-1pt}|}

\newcounter{defin}  \newcounter{lemma}  \newcounter{theorem}
\newcounter{property} \newcounter{corol}  \newcounter{remark} \newcounter{example}

\newenvironment{lemma}{\par\refstepcounter{lemma}%\noindent
     \textbf{Lemma \thelemma.} }{\rm\par}
\newenvironment{theorem}{\par\refstepcounter{theorem}%\noindent
     \textbf{Theorem \thetheorem.}\ }{\rm\par}
\newenvironment{property}{\par\refstepcounter{property}%\noindent
     \textbf{Proposition \theproperty.}\ }{\rm\par}
\newenvironment{corollary}{\par\refstepcounter{corol}%\noindent
     \textbf{Corollary \thecorol.} }{\rm\par}
\newenvironment{definition}{\par\refstepcounter{defin}%\noindent
     \textbf{Definition \thedefin.}\ }{\rm\par}
\newenvironment{remark}{\par\refstepcounter{remark}%\noindent
     \textbf{Remark \theremark.}}{\rm\par}

\begin{document}

\title{On completion of the cone of CP linear maps with respect to the energy-constrained diamond norm}
\author{M.E. Shirokov\footnote{Steklov Mathematical Institute, RAS, Moscow, email:msh@mi.ras.ru}}
\date{}
\maketitle

\vspace{-10pt}

\begin{abstract}
For a given positive operator  $G$ we consider the cones of linear maps  between Banach spaces of trace class operators characterized by the Stinespring-like representation with $\sqrt{G}$-bounded and $\sqrt{G}$-infinitesimal operators correspondingly. We prove the completeness of both cones w.r.t. the energy-constrained diamond norm induced by $G$ (as an energy observable) and the coincidence of the second cone with the completion of the cone of CP linear maps w.r.t. this norm.

We show that the sets of quantum channels and  quantum operations are complete w.r.t. the energy-constrained diamond norm for any  energy observable.

Some properties of the maps belonging to the introduced cones are described. In particular, the corresponding generalization of the Kretschmann-Schlingemann-Werner theorem is obtained.

We also give a nonconstructive description of the completion of the set of all Hermitian-preserving completely bounded linear maps w.r.t. the ECD norm.
\end{abstract}

\tableofcontents

\section{Introduction}

The norm of complete boundedness (typically called the diamond norm) on the set of completely positive (CP) linear maps between Banach spaces of trace class operators is widely used in the quantum
theory \cite{Kit,Watrous,Wilde}. The corresponding distance between quantum channels can be treated  as a  measure of distinguishability between these channels by quantum measurements \cite[Ch.9]{Wilde}. Nevertheless, the topology (convergence) generated  by  the diamond norm distance is, in general, too strong for description of physical perturbations of infinite-dimensional quantum channels \cite{SCT,W-EBN}. This is explained, briefly speaking, by the fact that
the definition of the diamond norm does not separate input states with finite energy (which can be produced in a physical experiment) and unrealizable states with infinite energy. To take  this "energy discrimination" of quantum states into account the energy-constrained diamond norms (ECD norms)
 are introduced independently in \cite{SCT,W-EBN}, where it is shown that these norms induce appropriate metric and topology on the set of quantum channels and operations (slightly different energy-constrained diamond norm is  used in \cite{Pir}). In particular, the ECD norms induce the strong convergence topology on the set of quantum channels and operations provided that the input system Hamiltonian is of discrete type \cite{SCT}.

The ECD norms turned out to be a useful tool for quantitative continuity analysis of basic capacities of energy-constrained infinite-dimensional channels \cite{SCT,W-EBN}. These norms are also used in study of quantum dynamical semigroups \cite{Datta,W-EBN,QDS}.

The cone $\F(A,B)$ of CP linear maps between Banach spaces $\T(\H_A)$ and $\T(\H_B)$ of trace class operators is complete w.r.t. the diamond norm metric, but it is  not complete w.r.t. the ECD norm metric induced by an unbounded Hamiltonian $G$ of the system $A$.
The aim of this paper is to describe the  completion of the cone $\F(A,B)$  w.r.t. the ECD norm metric induced by any positive operator $G$ on $\H_A$. We introduce the cone $\F_G(A,B)$  of linear maps defined on the linear span of states with finite energy and characterized by the Stinespring-like representation with $\sqrt{G}$-bounded operators from $\H_A$ to $\H_{B}\otimes\H_{E}$ (where $E$ is an environment). We prove that this cone is complete w.r.t. to the ECD norm metric and that
the  completion of the cone $\F(A,B)$  w.r.t. this metric coincides with the proper subcone $\F^0_G(A,B)$ of $\F_G(A,B)$ consisting of linear maps having Stinespring-like representation with $\sqrt{G}$-infinitesimal operators from   $\H_A$ to $\H_{B}\otimes\H_{E}$.
We also prove that the diamond norm bounded subset of $\F(A,B)$ is complete
w.r.t the ECD norm metric  if and only if it is closed w.r.t this metric. This implies completeness
of the sets of  quantum channels and quantum operations w.r.t. the ECD norm metric.

We obtain characterizations of the cones $\F_G(A,B)$ and $\F^0_G(A,B)$ in terms of the Kraus representation and generalize the Kretschmann-Schlingemann-Werner theorem  (obtained in \cite{Kr&W}) to all the maps from the cone $\F_G(A,B)$.

In the last section we give a nonconstructive description of the  completion $\Y_G(A,B)$ of the real linear space $\Y(A,B)$ of all Hermitian-preserving completely bounded linear maps from $\T(\H_A)$ to $\T(\H_B)$ w.r.t. the ECD norm. We prove that the positive cone in $\Y_G(A,B)$ coincides with the cone $\F^0_G(A,B)$ provided that the operator $G$ has a discrete spectrum of finite multiplicity.

\section{Preliminaries}

\subsection{Basic notations}

Let $\H,\H'$ be  separable infinite-dimensional Hilbert spaces, $\mathfrak{B}(\H,\H')$
-- the Banach space of all bounded operators from $\H$ to $\H'$ with the operator norm $\|\!\cdot\!\|$ and $\mathfrak{T}(\H,\H')$ --
the Banach space of all trace-class operators from $\H$ to $\H'$ with the trace norm $\|\!\cdot\!\|_1$ (the Schatten class of order 1) \cite{B&R,R&S}.
We will assume that $\mathfrak{B}(\H)\doteq\mathfrak{B}(\H,\H)$ and $\mathfrak{T}(\H)\doteq\mathfrak{T}(\H,\H)$.
Let $\mathfrak{T}_{+}(\H)$ be the cone of positive operators in
$\mathfrak{T}(\mathcal{H})$ and $\S(\H)$ its closed convex subset consisting of operators with unit trace (called \emph{quantum states}). Extreme points of $\S(\H)$ are 1-rank projectors called \emph{pure states}.

Trace-class operators will be usually denoted by the Greek letters $\rho $, $\sigma $, $\omega $, ... in contrast to other linear operators (bounded and unbounded) denoted by the Latin letters. The Greek letters $\phi$, $\varphi$, $\psi$,...  will be used for vectors in a Hilbert space.

Denote by $I_{\H}$ the unit operator on a Hilbert space
$\mathcal{H}$ and by $\id_{\mathcal{\H}}$ the identity
transformation of the Banach space $\mathfrak{T}(\mathcal{H})$.

The \emph{Bures distance} between operators $\rho$ and $\sigma$ in $\T_+(\H)$ is defined as
\begin{equation}\label{B-d-s}
  \beta(\rho,\sigma)=\sqrt{\|\rho\|_1+\|\sigma\|_1-2\sqrt{F(\rho,\sigma)}},
\end{equation}
where
\begin{equation}\label{F-d-s}
  F(\rho,\sigma)=\|\sqrt{\rho}\sqrt{\sigma}\|^2_1
\end{equation}
is the \emph{fidelity} of $\rho$ and $\sigma$.  The following relations  between the Bures distance
and the  trace-norm distance hold (cf. \cite{H-SCI,Wilde})
\begin{equation}\label{B-d-s-r}
\frac{\|\rho-\sigma\|_1}{\sqrt{\|\rho\|_1}+\sqrt{\|\sigma\|_1}}\leq\beta(\rho,\sigma)\leq\sqrt{\|\rho-\sigma\|_1}.
\end{equation}

If quantum systems $A$ and $B$ are described by Hilbert spaces  $\H_A$ and $\H_B$ then the bipartite system $AB$ is described by the tensor product of these spaces, i.e. $\H_{AB}\doteq\H_A\otimes\H_B$. A state in $\S(\H_{AB})$ is denoted by $\rho_{AB}$, its marginal states $\Tr_{\H_B}\rho_{AB}$ and $\Tr_{\H_A}\rho_{AB}$ are denoted, respectively, by $\rho_{A}$ and $\rho_{B}$ (here $\Tr_X\rho_{AB}\doteq\Tr_{\H_X}\rho_{AB}$).
\smallskip

We will pay a special attention to the class of  unbounded densely defined positive operators on $\H$ having discrete spectrum of finite multiplicity.
In Dirac's notations any such operator ${G}$ can be represented as follows
\begin{equation}\label{H-rep}
G=\sum_{k=0}^{+\infty} E_k|\tau_k\rangle\langle \tau_k|
\end{equation}
on the domain $\mathcal{D}(G)=\{ \varphi\in\H\,|\,\sum_{k=0}^{+\infty} E^2_k|\langle\tau_k|\varphi\rangle|^2<+\infty\}$, where
$\left\{\tau_k\right\}_{k=0}^{+\infty}$ is the orthonormal basis of eigenvectors of ${G}$
corresponding to the nondecreasing sequence $\left\{\smash{E_k}\right\}_{k=0}^{+\infty}$ of eigenvalues
tending to $+\infty$. We will use the following (cf.\cite{W-EBN})\smallskip

\begin{definition}\label{D-H}
An operator ${G}$ having representation (\ref{H-rep}) is called \emph{discrete}.
\end{definition}\smallskip

\subsection{Relatively bounded operators and the operator \emph{E}-norms}

In this paper we will consider  linear operators between separable Hilbert spaces $\H$ and $\H'$ relatively bounded with
respect to the  operator $\sqrt{G}$,  where $G$ is a positive semidefinite operator on $\H$ with  dense domain $\mathcal{D}({G})$ treated as a Hamiltonian (energy observable) of a quantum system associated with the space $\H$. We will assume that
\begin{equation}\label{G-cond}
\inf\left\{\shs\|G\varphi\|\,|\,\varphi\in\mathcal{D}({G}),\|\varphi\|=1\shs\right\}=0.
\end{equation}

A linear operator $A:\H\rightarrow\H'$ is called relatively bounded w.r.t. the operator $\sqrt{G}$ (briefly, $\sqrt{G}$-bounded) if
$\mathcal{D}(\sqrt{G})\subseteq \D(A)$ and
\begin{equation}\label{rb-rel}
\|A\varphi\|^2\leq a^2\|\varphi\|^2+b^2\|\sqrt{G}\varphi\|^2,\quad \forall \varphi\in\mathcal{D}(\sqrt{G}),
\end{equation}
for some  nonnegative numbers $a$ and $b$ \cite{Kato}. The $\sqrt{G}$-bound of $A$ (denoted by $b_{\sqrt{G}}(A)$ in what follows) is defined as
the infimum of the values $b$ for which (\ref{rb-rel}) holds with some $a$. If the $\sqrt{G}$-bound is equal to zero then $A$ is called $\sqrt{G}$-infinitesimal operator (infinitesimally bounded w.r.t. $\sqrt{G}$) \cite{Kato,R&S+,BS}.

Since $\sqrt{G}$ is a closed operator, the linear space $\D(\sqrt{G})$ equipped with the inner product
\begin{equation}\label{in-p}
\langle\varphi|\psi\rangle^G_E=\langle\varphi|\psi\rangle+\langle\varphi|G|\psi\rangle/E,\quad E>0,
\end{equation}
is a Hilbert space denoted by $\H_E^G$ in what follows \cite{R&S}. A restriction of any $\sqrt{G}$-bounded operator from $\H$ to $\H'$  to the set $\D(\sqrt{G})$ can be treated as a bounded operator from $\H_E^G$ to $\H'$ and, vise versa, any bounded operator from $\H_E^G$ to $\H'$ induces a $\sqrt{G}$-bounded operator from $\H$ to $\H'$. Thus, the linear space of all $\sqrt{G}$-bounded operators from $\H$ to $\H'$ equipped with the norm
\begin{equation}\label{eq-norms-2}
\sn A\sn^{G}_{E}=\sup_{ \varphi\in \D(\sqrt{G})}\frac{\|A\varphi\|}{\sqrt{\|\varphi\|^2+\|\sqrt{G}\varphi\shs\|^2/E}},\quad E>0,
\end{equation}
is a Banach space.\footnote{We identify operators coinciding on the set $\D(\sqrt{G})$.} For our purposes it is more convenient to  use the equivalent norm
\begin{equation}\label{ec-on}
 \|A\|^{G}_E\doteq \sup_{\substack{\rho\in\mathfrak{S}(\mathcal{H}):
\Tr G\rho \leq E}}\sqrt{\Tr A\rho A^*},\quad E>0,
\end{equation}
on the linear space of all $\sqrt{G}$-bounded operators from $\H$ to $\H'$. The supremum here is over all states $\rho$ in $\S(\H)$ such that  $\Tr G\rho\leq E$.\footnote{The value of $\Tr {G}\rho$ (finite or infinite) is defined as $\sup_n\Tr P_nG\rho$, where $P_n$ is the spectral projector of $G$ corresponding to the interval $[0,n]$.} By Lemma 5 in \cite{ECN}\footnote{In \cite{ECN}  linear operators on a single Hilbert space are considered. But all the results obtained therein are directly generalized to linear operators between different Hilbert spaces \cite[Remark 1]{ECN}.} for any $\sqrt{G}$-bounded operator $A:\H\rightarrow\H'$ the function $\rho\mapsto A\rho A^*$ is well defined on the set $\,\T^+_{G}(\H)\doteq \{\rho\in\T_{+}(\H)\,|\,\Tr {G}\rho<+\infty\shs\}\,$ by the expression
\begin{equation}\label{aa-d}
  A\rho A^*\doteq\sum_i|\alpha_i\rangle\langle\alpha_i|,\qquad |\alpha_i\rangle=A|\varphi_i\rangle,
\end{equation}
where $\rho=\sum_i |\varphi_i\rangle\langle\varphi_i|$ is any decomposition of $\rho\in\T^+_{G}(\H)$ into $1$-rank positive operators. This function is affine
and takes values in $\T_{+}(\H')$. So, the r.h.s. of (\ref{ec-on}) is well defined for any $\sqrt{G}$-bounded operator $A$. Due to condition (\ref{G-cond})
the supremum in (\ref{ec-on}) can be taken over all operators $\rho\in\T_+(\H)$ such that $\Tr G\rho\leq E$ and $\Tr \rho\leq 1$ \cite[Pr.3]{ECN}.\smallskip

The norm $\|\cdot\|^{G}_{E}$ called the operator \emph{E}-norm  in \cite{ECN} can be also defined by the following equivalent expressions
\begin{equation}\label{en-def-r-1}
 \!\!\textstyle \|A\|^{G}_E=\sup\left\{\sqrt{\sum_i\|A\varphi_i\|^2}\,\left|\,\{\varphi_i\}\subset\D(\sqrt{G}):\sum_i\|\varphi_i\|^2\leq 1,\;\sum_i\|\sqrt{G}\varphi_i\|^2\leq E\right.\right\}
\end{equation}
and
\begin{equation}\label{en-def-r-2}
\textstyle \|A\|^{G}_E=\sup\left\{\|A\otimes I_{\K}\varphi\|\,\left|\,\varphi\in\D(\sqrt{G}\otimes I_{\K}):\|\varphi\|\leq 1,\;\|\sqrt{G}\otimes I_{\K}\varphi\|^2\leq E\right.\right\},
\end{equation}
where $\K$ is a separable infinite-dimensional Hilbert space. If $G$ is a discrete  operator (Def.\ref{D-H}) then all the above expressions are simplified as follows
\begin{equation}\label{en-def-r-1++}
 \textstyle \|A\|^{G}_E=\sup\left\{\|A\varphi\|\,\left|\,\varphi\in\D(\sqrt{G}):\|\varphi\|\leq 1,\;\|\sqrt{G}\varphi\|^2\leq E\right.\right\}.
\end{equation}
Validity of the simplified expression (\ref{en-def-r-1++}) in the case of arbitrary positive operator $G$ is an interesting open question (see the Appendix in \cite{ECN}).\footnote{The r.h.s. of (\ref{en-def-r-1++}) defines  a norm on the set of all $\sqrt{G}$-bounded operators denoted by $\|\cdot\|^{G}_{\circ,E}$ in \cite{ECN}, where it is shown that this norm is equivalent to the norms $\|\cdot\|^{G}_{E}$ and $\sn \cdot\sn^{G}_{E}$ and that the function $E\mapsto\left[\|A\|^{G}_{E}\right]^2$ is the concave envelope (hull) of the function $E\mapsto\left[\|A\|^{G}_{\circ,E}\right]^2$. It means, in particular, that the conjectured coincidence of $\|\cdot\|^{G}_{E}$ and $\|\cdot\|^{G}_{\circ,E}$ is equivalent to concavity of the function   $E\mapsto\left[\|A\|^{G}_{\circ,E}\right]^2$ for any $\sqrt{G}$-bounded operator $A$.}

For any $\sqrt{G}$-bounded operator $A:\H\rightarrow\H'$ both norms $\|A\|^{G}_E$ and $\sn A\sn^{G}_E$ are nondecreasing functions of $E$ tending to $\|A\|\leq+\infty$ as $E\rightarrow+\infty$. For any $E>0$ they are related by the inequalities
\begin{equation}\label{eq-one}
\sqrt{1/2}\|A\|^{G}_{E}\leq \sn A\sn^{G}_{E}\leq\|A\|^{G}_{E},
\end{equation}
which show the equivalence of these norms on the
set of all $\sqrt{G}$-bounded operators \cite{ECN}. Moreover, for any $\sqrt{G}$-bounded operator $A$ the functions
$E\mapsto \|A\|^{G}_{E}$ and $E\mapsto \sn A\sn^{G}_{E}$ are completely determined by each other via the following expressions (\cite[Th.3A]{ECN}):
\begin{equation}\label{s-tr}
\sn A\sn^{G}_{E}=\sup_{t>0}\|A\|^{G}_{tE}/\sqrt{1+t},\quad  \|A\|^{G}_{E}=\inf_{t>0}\sn A \sn^{G}_{tE}\sqrt{1+1/t},\quad E>0.
\end{equation}

One of the main advantages of the norm $\|\cdot\|^{G}_E$ is the concavity of the function $E\mapsto\left[\|A\|^{G}_{E}\right]^p$ for any $p\in(0,2]$ and any $\sqrt{G}$-bounded operator $A$ which essentially simplifies quantitative analysis of functions depending on $\sqrt{G}$-bounded operators \cite[Section 5]{ECN}.\footnote{The function $E\mapsto\left[\sn A\sn^{G}_{E}\right]^p$ is not concave in general for any $p\in(0,2]$ \cite[Section 3.1]{ECN}.}
This property implies, in particular, that
\begin{equation}\label{E-n-eq}
\|A\|^{G}_{E_1}\leq \|A\|^{G}_{E_2}\leq \sqrt{E_2/E_1}\|A\|^{G}_{E_1}\quad\textrm{ for any } E_2>E_1>0.
\end{equation}
Hence for given operator ${G}$ all the norms $\|\!\cdot\!\|^{G}_{E}$, $E>0$, are equivalent on the
set of all $\sqrt{G}$-bounded operators. By inequalities (\ref{eq-one})
the same is true for the norms $\sn\!\cdot\!\sn^{G}_{E}$, $E>0$.

Another advantage of the norm $\|A\|^{G}_E$  essentially used in this paper is connected with its appearance in the generalized version of
Kretschmann-Schlingemann-Werner theorem  described in Section 4.3.

Denote by $\B_G(\H,\H')$ the linear space of all $\sqrt{G}$-bounded operators from $\H$ into $\H'$ equipped with any of the equivalent norms $\|\cdot\|^{G}_E$, $E>0$. The equivalence of the norms $\|\cdot\|^{G}_E$ and $\sn \cdot\sn^{G}_E$ mentioned before implies that $\B_G(\H,\H')$ is a (nonseparable) Banach space. The $\sqrt{G}$-bound $b_{\sqrt{G}}(\cdot)$ is a continuous seminorm on $\B_G(\H,\H')$, for any operator $A\in\B_G(\H,\H')$ it can be determined by the formula
\begin{equation}\label{G-bound}
b_{\sqrt{G}}(A)=\lim_{E\rightarrow+\infty}\|A\|^{G}_{E}/\sqrt{E},
\end{equation}
where the limit can be replaced by infimum over all $E>0$ \cite[Theorem 3B]{ECN}.\smallskip

The closed subspace $\,\B^0_{\!G}(\H,\H')$ of $\,\B_{\!G}(\H,\H')$ consisting of all $\sqrt{G}$-infinitesimal operators, i.e. operators with
the $\sqrt{G}$-bound equal to $0$, coincides with the completion of $\B(\H,\H')$ w.r.t. any of the norms $\|\!\cdot\!\|^{G}_E$, $E>0$ \cite[Theorem 3C]{ECN}. It follows from (\ref{G-bound}) that an operator $A$ belongs to the space $\B^0_{\!G}(\H,\H')$ if and only if
\begin{equation}\label{s-cond}
 \|A\|^{G}_E=o\shs(\sqrt{E})\quad\textup{ as }\quad E\rightarrow+\infty.
\end{equation}

In the following two lemmas we collect some results from Section 5 in  \cite{ECN} which will be used in this paper.\smallskip

\begin{lemma}\label{vbl} \emph{If $A$ is a $\sqrt{G}$-bounded operator from $\H$ to $\H'$ then for any separable Hilbert space $\K$  the operator $A\otimes I_{\K}$ naturally defined on the set $\,\D(\sqrt{G})\otimes\K$  has a unique linear $\sqrt{G}\otimes I_{\K}$-bounded extension to the set $\,\D(\sqrt{G}\otimes I_{\K})$.\footnote{$\D(\sqrt{G})\otimes \K$ is the linear span of all the vectors $\varphi\otimes\psi$, where $\varphi\in\D(\sqrt{G})$ and $\psi\in\K$.} This extension (also denoted by $A\otimes I_{\K}$) has the following property
\begin{equation}\label{s-prop}
 A\otimes I_{\K}\!\left(\sum_{i}|\varphi_i\rangle\otimes|\psi_i\rangle\right)=\sum_{i}A|\varphi_i\rangle\otimes|\psi_i\rangle
\end{equation}
for any countable sets $\{\varphi_i\}\subset\D(\sqrt{G})$ and $\{\psi_i\}\subset\K$ such that $\sum_{i}\|\sqrt{G}\varphi_i\|^2<+\infty$ and $\langle\psi_i|\psi_j\rangle=\delta_{ij}$, which implies that $\|A\otimes I_{\K}\|^{{G}\otimes I_{\K}}_E=\|A\|^{{G}}_E$ for any $E>0$.}\smallskip

\emph{If $A$ is a $\sqrt{G}$-infinitesimal operator from $\H$ to $\H'$ then the extension of $A\otimes I_{\K}$ mentioned above is uniformly  continuous on the set
\begin{equation}\label{s-3}
\V_{E}\doteq \{\shs\eta\in\D(\sqrt{G}\otimes I_{\K})\,|\, \|\sqrt{G}\otimes I_{\K}\eta\|^2\leq E\shs\}
\end{equation}
for any $E>0$. Quantitatively,
\begin{equation*}%\label{v-CB}
\|A\otimes I_{\K}(\eta-\theta)\|\leq f_A(E,\varepsilon)
\end{equation*}
for any vectors $\,\eta$ and $\,\theta$ in $\V_{E}$ such that $\|\eta-\theta\|\leq\varepsilon$, where
$f_A(E, \varepsilon)=\varepsilon\|A\|^{G}_{4E/\varepsilon^2}$ is a function vanishing as $\,\varepsilon\rightarrow 0^+$ by condition (\ref{s-cond}).}%\smallskip

%\emph{If $\,A\in\B_{\!G}(\H)\setminus\B^0_{\!G}(\H)$ then the operator $A\otimes I_{\K}$ is not continuous  on the set $\,\V_{E}$ for any $E>0$.}
\end{lemma}\medskip

%\begin{remark}\label{w-c} Property (\ref{s-prop}) implies that
%$
%(A\otimes I_{\K})(I_{\H}\otimes W)|\varphi\rangle=(I_{\H}\otimes W)(A\otimes I_{\K})|\varphi\rangle
%$
%for any $\varphi\in\D(\sqrt{G}\otimes I_{\K})$ and a  partial isometry $W\in\B(\K)$ s.t.
%$I_{\H}\otimes W^{*}W |\varphi\rangle=|\varphi\rangle$.
%\end{remark}\medskip

\begin{lemma}\label{qsl} \emph{For any $\sqrt{G}$-bounded operators $A$ and $B$ from $\H$ to $\H'$
the  function $\shs\rho\mapsto A\rho B^*\in\T(\H')$ is well defined on the set $\,\T^+_{G}(\H)\doteq \{\rho\in\T_{+}(\H)\,|\,\Tr {G}\rho<+\infty\shs\}\,$ by the formula \footnote{We define the operator $A\rho B^*$ in such a way to avoid the notion of adjoint operator, since we make no assumptions about closability of the operators $A$ and $B$.}
\begin{equation}\label{ab-d}
  A\rho B^*\doteq\sum_i|\alpha_i\rangle\langle\beta_i|,\qquad |\alpha_i\rangle=A|\varphi_i\rangle,\;|\beta_i\rangle=B|\varphi_i\rangle,
\end{equation}
where $\rho=\sum_i |\varphi_i\rangle\langle\varphi_i|$ is any decomposition of $\,\rho\in\T^+_{G}(\H)$ into $1$-rank positive operators}

\noindent\emph{For any operator $\rho$ in $\,\T^+_{G}(\H)$ such that $\,\Tr\rho\leq 1$ the following inequality holds}
$$
\|A\rho B^*\|_1 \leq \|A\|^{G}_{E_{\rho}}\|B\|^{G}_{E_{\rho}},\quad\textit{where}\quad E_{\rho}=\Tr G\rho.
$$
\emph{For any $\sqrt{G}$-infinitesimal operators $A$ and $B$ the function $\shs\rho\mapsto A\rho B^*$  is uniformly  continuous on the set
$\,\mathfrak{C}_{{G},E}\doteq \{\rho\in\T_{+}(\H)\,|\,\Tr\rho\leq 1,  \Tr {G}\rho\leq E\shs\}$ for any $E>0$. Quantitatively,
\begin{equation*}%\label{ab-cb}
\|A\rho B^*- A\sigma B^*\|_1\leq \|A\|_E^G f_B(E,\sqrt{\varepsilon})+\|B\|_E^G f_A(E,\sqrt{\varepsilon})
\end{equation*}
for any operators $\,\rho$ and $\sigma$ in $\,\C_{{G},E}$ such that $\|\rho-\sigma\|_1\leq\varepsilon$, where
$f_X$ is the function defined in Lemma  \ref{vbl}.}
\end{lemma}\smallskip

The first assertion of Lemma \ref{qsl} implies that $\shs\rho\mapsto A\rho B^*$ is an affine function on the cone $\,\T^+_{G}(\H)$.

\section{The main results}

For a completely positive (CP)  linear map $\,\Phi:\T(\H_A)\rightarrow \T(\H_B)\,$ the Stinespring theorem (cf.\cite{St}) implies existence of a separable Hilbert space
$\mathcal{H}_E$ and a bounded operator
$V_{\Phi}:\mathcal{H}_A\rightarrow\mathcal{H}_{BE}\doteq\mathcal{H}_B\otimes\mathcal{H}_E$ such
that
\begin{equation}\label{St-rep}
\Phi(\rho)=\mathrm{Tr}_{E}V_{\Phi}\rho V_{\Phi}^{*},\quad
\rho\in\mathfrak{T}(\mathcal{H}_A),
\end{equation}
where $\Tr_E$ denotes the partial trace over $\H_E$. If $\Phi$ is trace-preserving (correspondingly, trace-non-increasing) then $V_{\Phi}$ is an isometry
(correspondingly, contraction) \cite[Ch.6]{H-SCI}.

The norm of complete boundedness of a normal linear map between the algebras $\B(\H_B)$ and $\B(\H_A)$ (cf. \cite{Paul}) induces (by duality) the
diamond norm
\begin{equation}\label{d-n-def}
\|\Phi\|_{\diamond}\doteq\sup_{\omega\in\T(\H_{AR}),\|\omega\|_1\leq 1}\|\Phi\otimes \id_R(\omega)\|_1
\end{equation}
on the set of all linear maps between the Banach spaces  $\T(\H_A)$ and $\T(\H_B)$, where $\H_R$ is a separable Hilbert space and $\H_{AR}=\H_{A}\otimes \H_{R}$ \cite{Kit}. If $\Phi$ is a Hermitian preserving map
then the supremum in (\ref{d-n-def}) can be taken over the set $\S(\H_{AR})$ \cite[Ch.3]{Watrous}.
\smallskip

The diamond norm is widely used in the quantum  theory, but in general the convergence induced by this norm is too strong for description of
physical perturbations of infinite-dimensional quantum channels: there exist quantum channels with close physical parameters such that the diamond norm distance between them is equal to $2$ \cite{W-EBN}. The reason of this inconsistency is pointed briefly in the Introduction. By taking it into account the energy-constrained  diamond norms (called ECD norms in what follows)
\begin{equation}\label{E-sn}
  \|\Phi\|_{\diamond,E}^G\doteq\sup_{\omega\in\S(\H_{AR}):\Tr {G}\omega_A\leq E}\|\Phi\otimes \id_R(\omega)\|_1,\quad E>0,
\end{equation}
on the set  $\Y(A,B)$ of Hermitian-preserving
linear maps from $\T(\H_{A})$ to $\T(\H_{B})$ are introduced independently in \cite{SCT} and \cite{W-EBN} (here  $G$ is a positive operator on the space $\H_A$ satisfying condition (\ref{G-cond}) treated as a Hamiltonian of a quantum system $A$).\footnote{Slightly different energy-constrained diamond norm is used in \cite{Pir}.}
\smallskip

In \cite{W-EBN} it is shown that for any given $\Phi\in\Y(A,B)$ the nondecreasing nonnegative function $E\mapsto\|\Phi\|_{\diamond,E}^G$ is concave on $\mathbb{R}_+$. This implies that
\begin{equation}\label{ECD-n-eq}
\|\Phi\|_{\diamond,E_1}^G\leq \|\Phi\|_{\diamond,E_2}^G\leq (E_2/E_1)\|\Phi\|_{\diamond,E_1}^G\quad\textrm{ for any } E_2>E_1>0.
\end{equation}
Hence for given operator ${G}$ all the norms $\|\!\cdot\!\|_{\diamond,E}^G$, $E>0$, are equivalent on $\Y(A,B)$.

In \cite{SCT,W-EBN} it is shown that the ECD norms induce adequate metric on the set of infinite-dimensional quantum channels which is consistent with the energy separations of quantum states. In particular, this metric  has an operational interpretation in terms of discriminating quantum channels with test states of bounded energy \cite{W-EBN}. If ${G}$ is a discrete unbounded operator (see Def.\ref{D-H}) then  any of the norms (\ref{E-sn}) generates the strong  convergence on the set of quantum channels \cite{SCT}.\footnote{The strong convergence of a sequence $\{\Phi_n\}$ of channels  to a channel $\Phi_0$  means that
$\lim_{n\rightarrow\infty}\Phi_n(\rho)=\Phi_0(\rho)\,\textup{ for all }\rho\in\S(\H_A)$.} This holds, for example, if $G$ is the Hamiltonian of a multi-mode quantum oscillator \cite[Ch.12]{H-SCI}. \smallskip

Denote by $\F(A,B)$ the cone of all CP linear maps from $\T(\H_A)$ into $\T(\H_B)$. The cone $\F(A,B)$ is complete  w.r.t. the metric induced by the diamond norm but not complete w.r.t. the metric induced by the ECD norm (provided that $G$ is an unbounded operator).
Our aim is to describe the completion of  the cone $\F(A, B)$ and of its important subsets
w.r.t. the metric induced by the ECD norm.\smallskip

Let $\H_E$ be a separable Hilbert space and $V$ an arbitrary $\sqrt{G}$-bounded operator from $\H_A$ to $\H_{BE}$. Lemma \ref{qsl} in Section 2.2 implies that $\rho\mapsto V\rho V^*$ is an affine map from the set $\S_G(\H_A)$ of all states $\rho$ in $\S(\H_A)$ with finite energy $\Tr G\rho$ into the cone $\T_+(\H_{BE})$ correctly defined by the formula
\begin{equation}\label{aa-d}
  V\rho V^*\doteq\sum_i|V\varphi_i\rangle\langle V \varphi_i|,
\end{equation}
where $\rho=\sum_i |\varphi_i\rangle\langle\varphi_i|$ is any decomposition of $\rho$ into 1-rank positive operators.
So, we may define the affine map
\begin{equation}\label{map}
\Phi(\rho)=\Tr_E V\rho V^*
\end{equation}
from the set $\S_G(\H_A)$ into the cone $\T_+(\H_{B})$ which can be extended to a unique linear map from
the linear span  $\T_G(\H_A)$ of $\S_G(\H_A)$ into the space $\T(\H_{B})$.\smallskip

We will say that two maps having form (\ref{map}) are $G$-\emph{equivalent} if they coincides on the set $\S_G(\H_A)$.
In what follows we will identify $G$-equivalent maps.\smallskip

Let $\F_G(A,B)$ be the set of all maps having form (\ref{map}) for some  separable Hilbert space $\H_E$ and operator $V\in\B_G(\H_A,\H_{BE})$.\footnote{$\B_G(\H_A,\H_{BE})$ is the Banach space of $\sqrt{G}$-bounded operators from $\H_A$ to $\H_{BE}$ (see Section 2.2).} By the Stinespring representation (\ref{St-rep}) the cone $\F(A,B)$ is naturally embedded into the set $\F_G(A,B)$ (the map $\Phi$ defined in (\ref{map}) belongs to the cone $\F(A,B)$ if and only if  the operator $V$ belongs to the subset $\B(\H_A,\H_{BE})$ of $\B_G(\H_A,\H_{BE})$).\smallskip

\begin{definition}\label{ro-def} The operator $V$ in any representation (\ref{map}) of a map $\Phi\in\F_G(A,B)$ will be called \emph{representing operator} for this map.
\end{definition}\smallskip

Let $\Phi$ be any map in $\F_G(A,B)$ with representing operator $V\in\B_G(\H_A,\H_{BE})$. By Lemma \ref{vbl} for any separable Hilbert space $\H_R$  the operator $V\otimes I_R$ belongs to the space $\B_{G\otimes I_R}(\H_{AR},\H_{BER})$ and $\|V\otimes I_R\|_E^{G\otimes I_R}=\|V\|_E^{G}$.
By Lemma \ref{qsl} the map
\begin{equation}\label{p-map}
\Theta(\omega)=\Tr_E [V\otimes I_R]\,\omega\, [V\otimes I_R]^*
\end{equation}
is well defined on the set
$\S_{G\otimes I_R}(\H_{AR})$ by the formula similar to (\ref{aa-d}). This map
does not depend on the representing operator $V$. Indeed, let $V'\in\B_G(\H_A,\H_{BE'})$ be another representing operator for $\Phi$
and $\Theta'$ the corresponding map (\ref{p-map}). It suffices to show that $\Theta$ and $\Theta'$ coincide at a pure state $|\eta\rangle\langle\eta|$, where $\eta$ is a any unit vector in $\D(\sqrt{G}\otimes I_R)$ (since any state in $\S_{G\otimes I_R}(\H_{AR})$
is decomposed into a convex combination of such pure states).

Any unit vector $\eta$  in $\D(\sqrt{G}\otimes I_R)$ has the representation
$$
|\eta\rangle=\sum_{i}|\varphi_i\rangle\otimes|\psi_i\rangle,
$$
where $\{\varphi_i\}$ and $\{\psi_i\}$ are sets  of vectors in $\D(\sqrt{G})$ and $\H_R$  correspondingly such that $\sum_{i}\|\sqrt{G}\varphi_i\|^2<+\infty$ and $\langle\psi_i|\psi_j\rangle=\delta_{ij}$. By Lemma \ref{vbl} in Section 2.2 we have
$$
V\otimes I_{R}|\eta\rangle=\sum_{i}V|\varphi_i\rangle\otimes|\psi_i\rangle\quad \textrm{and} \quad
 V'\otimes I_{R}|\eta\rangle=\sum_{i}V'|\varphi_i\rangle\otimes|\psi_i\rangle.
$$
Hence the continuity of a partial trace implies that
$$
\Theta(|\eta\rangle\langle\eta|)=\sum_{i,j}\Phi(|\varphi_i\rangle\langle\varphi_j|)\otimes|\psi_i\rangle\langle\psi_j|=\Theta'(|\eta\rangle\langle\eta|).
$$

Thus, for any map $\Phi$ in $\F_G(A,B)$ and  any separable Hilbert space $\H_R$ there is a unique map $\Phi\otimes \id_R$
in $\F_G(AR,BR)$ defined in (\ref{p-map}). This property can be treated as \emph{complete positivity} of $\Phi$. It implies that for any map $\Phi$ in the real linear span of $\F_G(A,B)$ we may define the ECD-norm $\|\Phi\|_{\diamond,E}^G$ by formula (\ref{E-sn}). So, for any two maps $\Phi$ and $\Psi$ in $\F_G(A,B)$ we may define the distance
\begin{equation}\label{ecd}
  D^G_{E}(\Phi,\Psi)\doteq \|\Phi-\Psi\|_{\diamond,E}^G,\quad E>0.
\end{equation}
It is clear  that $D_E^G$ is a metric on the set $\F_G(A,B)$ for any $E>0$ coinciding with the  ECD-norm
metric on the cone $\F(A,B)$.\smallskip

Note that $\F_G(A,B)$ is a cone as well. Indeed, if $\Phi\in\F_G(A,B)$ then the map $\lambda\Phi$ obviously
belongs to the set $\F_G(A,B)$ for any positive $\lambda$. If $\Phi_k(\rho)=\Tr_{E_k} V_k\rho V_k^*$, where $V_k\in\B_G(\H_{A},\H_{BE_k})$, $k=1,2$, then
$$
(\Phi_1+\Phi_2)(\rho)=\Tr_{\tilde{E}}\tilde{V}\rho \tilde{V}^*,\quad \rho\in\T_G(\H_A),
$$
where $\H_{\tilde{E}}=\H_{E_1}\oplus\H_{E_2}$ and $\tilde{V}$ is the operator from $\D(\sqrt{G})$ into $\H_{B\tilde{E}}=\H_{BE_1}\oplus\H_{BE_2}$ defined by setting $\tilde{V}|\varphi\rangle=V_1|\varphi\rangle\oplus V_2|\varphi\rangle$ for any $\varphi\in\D(\sqrt{G})$. It is easy to see that $\tilde{V}$ belongs to the set
$\B_G(\H_{A},\H_{B\tilde{E}})$. Thus, the map $\Phi_1+\Phi_2$ belongs to the set $\F_G(A,B)$.\smallskip

%The maps in the cone $\F_G(A,B)$ will be called \emph{$G$-bounded CP linear maps from $\T(\H_A)$ to $\T(\H_B)$}.\smallskip

Let $\,\F^0_G(A,B)$ be the subset of $\,\F_G(A,B)$ consisting of  maps having form (\ref{map}) for some  separable Hilbert space $\H_E$ and operator $V\in\B^0_G(\H_A,\H_{BE})$.\footnote{$\B_G^0(\H_A,\H_{BE})$ is the subspace of $\B_G(\H_A,\H_{BE})$ consisting of $\sqrt{G}$--infinitesimal operators, it coincides with the completion of the set $\B(\H_A,\H_{BE})$ w.r.t. the norm $\|\cdot\|_{E}^G$ (see Section 2.2).} The above arguments implies that $\F^0_G(A,B)$ is a subcone of $\F_G(A,B)$.
%The maps in $\F^0_G(A,B)$ will be called \emph{$G$-infinitesimal CP linear maps from $\T(\H_A)$ to $\T(\H_B)$}.
\smallskip

In the following theorem we assume that $G$ is a positive (semidefinite) operator on the space $\H_A$ satisfying condition (\ref{G-cond}). We
use the ECD norm metric $D_{E}^G$ defined in (\ref{ecd}) and the operator \emph{E}-norms $\|\cdot\|_{E}^G$ described in Section 2.2. \smallskip

\begin{theorem}\label{main}
A) \emph{The cone $\,\F_G(A,B)$ is complete w.r.t. the metric $D_E^G$.} \smallskip

\noindent B) \emph{The subcone $\,\F^0_G(A,B)$ coincides with the completion of the cone $\,\F(A,B)$ w.r.t. the metric $D_E^G$. Any  map $\,\Phi$ in $\F^0_G(A,B)$ with representing operator $\,V\in\B_G^0(\H_A,\H_{BE})$
is approximated in the metric $D_E^G$ by the sequence of maps
$$
\Phi_n(\rho)=\Tr_E VP_{n}\rho P_{n}V^*
$$
from the cone $\,\F(A,B)$ determined by a sequence $\{E_n\}\subset\mathbb{R}_+$ tending to $+\infty$, where $P_{n}$ is the spectral projector of $\,G$ corresponding to the interval $\,[0,E_n]$. Quantitatively,}\footnote{Since the operator $V$ satisfies condition (\ref{s-cond}), the r.h.s. of (\ref{c-rate}) tends to zero as $n\rightarrow+\infty$.}
\begin{equation}\label{c-rate}
D_E^G(\Phi_n,\Phi)\leq 2\sqrt{E/E_n}\|V\|_{E_n}^G\|V\|^G_E\quad\textit{ for all }\;n\; \textit{ such that }\; E_n\geq E.
\end{equation}

\emph{A map $\Phi$ in $\,\F_G(A,B)$ belongs to the subcone $\,\F^0_G(A,B)$ if and only if}
\begin{equation}\label{ph-cond}
\|\Phi\|_{\diamond,E}^G=o\shs(E)\quad \emph{as} \quad E\rightarrow+\infty.
\end{equation}

\noindent C) \emph{A diamond norm bounded subset of $\,\F(A,B)$ is complete w.r.t. the metric $D_E^G$ if and only if it is closed w.r.t. this metric.}
\end{theorem}\smallskip

The last assertion of Theorem \ref{main} implies the following \smallskip
\begin{corollary}\label{main-c}
\emph{The set of quantum channels and the set of quantum operations are complete w.r.t. the metric $D_E^G$ induced by any positive operator $G$.}
\end{corollary}\smallskip

If $G$ is a discrete operator (Def.\ref{D-H}) then the ECD norm generates the strong convergence topology
on the set of quantum channels and operations. So, in this case the assertion of Corollary \ref{main-c}
agrees with Theorem 1 in \cite{Ende}. \smallskip

\emph{Proof of Theorem \ref{main}.} A)  Assume that $\{\Phi_n\}$ is a Cauchy sequence in  $\F_G(A,B)$.
For a given $E>0$ take a subsequence $\{\Phi_{n_k}\}$ such that $D_E^G(\Phi_{n_k},\Phi_{n_{k+1}})\leq 4^{-k}$.
Lemma \ref{r-in+} below and the right inequality in (\ref{D-B-in}) imply existence of a separable Hilbert space $\H_{E}$ and a sequence  $\{V_k\}$ of operators
in $\B_G(\H_A,\H_{BE})$ such that $\Phi_{n_k}(\rho)=\Tr_{E}V_k\rho V_k^*$ for all $\rho\in\S_G(\H_A)$ and
$$
\|V_{k}-V_{k+1}\|^G_E\leq \sqrt{D_E^G(\Phi_{n_k},\Phi_{n_{k+1}})}\leq 2^{-k}\quad \forall k.
$$
Thus, $\{V_k\}$ is a Cauchy sequence in the Banach space $\B_G(\H_A,\H_{BE})$ and hence it has a limit $V_0\in\B_G(\H_A,\H_{BE})$. Lemma \ref{l-in} below  implies that the subsequence
$\{\Phi_{n_k}\}$  converges to the map $\Phi_0(\rho)\doteq\Tr_{E}V_0\rho V_0^*$ w.r.t. the metric $D_E^G$.
Hence the whole sequence $\{\Phi_{n}\}$ converges to the map $\Phi_0$ w.r.t. the metric $D_E^G$ as well.\smallskip

B) Since $\B^0_G(\H_A,\H_{BE})$ is a closed subspace of $\B_G(\H_A,\H_{BE})$, the arguments from the proof of part $A$ imply that
$\F^0_G(A,B)$ is a closed subcone of $\F_G(A,B)$.

To show the density of $\F(A,B)$ in $\F_G(A,B)$ it suffices to prove (\ref{c-rate}), since Lemma 6 in \cite{ECN} implies that all the operators $VP_{n}$ are bounded and hence all the maps $\Phi_n$ belong to the cone $\F(A,B)$.

For given $n$ let $\rho$ be  a  state in $\S(\H_A)$ such that $\Tr {G}\rho\leq E$  and $x_n=1-\Tr P_{n}\rho>0$. Let $\bar{P}_{n}=I_{A}-P_{n}$ and $\rho_n=x^{-1}_n\bar{P}_{n}\rho\bar{P}_{n}$. We have
$$
\Tr V\bar{P}_{n} \rho\bar{P}_{n} V^*=x_n\Tr V \rho_n V^*\leq x_n\left[\|V\|^{G}_{E/x_n}\right]^2\leq  (E/E_n)\left[\|V\|^{G}_{E_n}\right]^2.
$$
The first inequality follows from the definition of the norm $\|\cdot\|^{G}_E$ and the inequality
$\Tr {G}\rho_n\leq E/x_n$, the second one follows from concavity of the  function $E\mapsto \left[\|V\|^{G}_E\right]^2$, Lemma \ref{WL} below and the inequality $x_n\leq E/E_n$ (which holds, since $\Tr {G}\rho\leq E$). The above estimate implies that
$\|V-VP_{n}\|^{G}_{E}\leq \sqrt{E/E_n}\|V\|^{G}_{E_n}$.
So, inequality (\ref{c-rate}) follows from Lemma \ref{l-in} below (since it is easy to see that $\|VP_{n}\|^{G}_E\leq\|V\|^{G}_E$ for all $n$).

If $V$ is a representing operator for $\Phi$ then $\left[\|V\|_{E}^G\right]^2=\|\Phi\|_{\diamond,E}^G$ for all $E>0$. So,
the property (\ref{ph-cond}) characterising the maps in $\F^0_G(A,B)$ follows from  the property (\ref{s-cond}) characterising
the operators in $\B^0_G(\H_A,\H_{BE})$.
$\square$\smallskip

\begin{corollary}\label{main-c+}
\emph{If $\,\Phi\in\F^0_G(A,B)$ then any representing operator for $\Phi$ is $\sqrt{G}$-infini- tesimal, i.e. belongs to the space $\B^0_G(\H_A,\H_{BE})$ for some Hilbert space $\H_E$.}
\end{corollary}\smallskip

\emph{Proof.} Since $\left[\|V\|_{E}^G\right]^2=\|\Phi\|_{\diamond,E}^G$ for any representing operator $V$ for $\Phi$, this assertion follows from
the characterising properties (\ref{s-cond}) and (\ref{ph-cond}). $\square$ \smallskip

\begin{lemma}\label{WL} \cite{W-CB} \emph{If $f$ is a concave nonnegative function on $[0,+\infty)$ then for any positive $x< y$ and any $z\geq0$ the inequality $\,xf(z/x)\leq yf(z/y)\,$ holds.}
\end{lemma}\smallskip

\begin{lemma}\label{l-in} \emph{Let $\Phi$ and $\Psi$ be any maps in $\F_G(A,B)$ with representing operators  $V_{\Phi}$ and $V_{\Psi}$ in $\B_G(\H_A,\H_{BE})$. Then}
$$
D_E^G(\Phi,\Psi)\leq \|V_{\Phi}-V_{\Psi}\|^G_E\left(\|V_{\Phi}\|_E^G+\|V_{\Psi}\|_E^G\right).
$$
\end{lemma}

\emph{Proof.} Let $\omega$ be any state in $\S(\H_{AR})$ such that $\Tr {G}\omega_A\leq E$. We have\footnote{We write $V_{\Phi}\otimes I_R \cdot\omega\cdot [V_{\Phi}\otimes I_R]^*$ instead of $V_{\Phi}\otimes I_R \cdot\omega\cdot V^*_{\Phi}\otimes I_R$ to emphasise that this operator is defined by formula (\ref{ab-d}) not assuming existence of densely defined adjoint operator $V^*_{\Phi}$.}
$$
\begin{array}{c}
\|\Phi\otimes \id_R(\omega)-\Psi\otimes \id_R(\omega)\|_1\leq \|V_{\Phi}\otimes I_R \cdot\omega\cdot [V_{\Phi}\otimes I_R]^*-V_{\Psi}\otimes I_R \cdot\omega\cdot [V_{\Psi}\otimes I_R]^*\|_1\\\\
\leq\|(V_{\Phi}-V_{\Psi})\otimes I_R \cdot\omega\cdot [V_{\Phi}\otimes I_R]^*\|_1+\|V_{\Psi}\otimes I_R \cdot\omega\cdot[(V_{\Phi}-V_{\Psi})\otimes I_R]^*\|_1\\\\
\leq\|(V_{\Phi}-V_{\Psi})\otimes I_R\|^{G\otimes I_R}_E \|V_{\Phi}\otimes I_R\|^{G\otimes I_R}_E+\|(V_{\Phi}-V_{\Psi})\otimes I_R\|^{G\otimes I_R}_E\|V_{\Psi}\otimes I_R\|^{G\otimes I_R}_E
\\\\
=\|V_{\Phi}-V_{\Psi}\|^G_E\, \|V_{\Phi}\|^G_E+\|V_{\Phi}-V_{\Psi}\|^G_E\, \|V_{\Psi}\|^G_E.
\end{array}
$$
The first and the second inequalities follow from the properties of the trace norm (the non-increasing under partial trace and the triangle inequality), the third inequality follows from Lemma \ref{qsl}, the last equality  -- from Lemma \ref{vbl}. $\square$\smallskip

The Bures distance between any maps $\Phi$  and $\Psi$ in $\F(A,B)$ is defined as
\begin{equation}\label{b-dist}
\beta(\Phi,\Psi)=\sup_{\omega\in\S(\H_{AR})}\beta(\Phi\otimes \id_R(\omega),\Psi\otimes \id_R(\omega)),
\end{equation}
where $\beta$ in the r.h.s. is the Bures distance between positive trace class operators defined in (\ref{B-d-s}) and $\H_R$ is an infinite-dimensional separable Hilbert space \cite{B&Co,Kr&W}. \smallskip

Following \cite{CID} consider  the \emph{energy-constrained Bures distance}
\begin{equation}\label{ec-b-dist}
\beta_E^G(\Phi,\Psi)=\sup_{\omega\in\S(\H_{AR}), \Tr H_A\omega_A\leq E} \beta(\Phi\otimes \id_R(\omega),\Psi\otimes \id_R(\omega)), \quad E>0,
\end{equation}
between any maps $\Phi$  and $\Psi$ in $\F(A,B)$.\footnote{Estimations of this distance for real quantum channels can be found in \cite{Nair}.} The arguments before (\ref{ecd}) show that this distance is well defined by the same formula for any maps $\Phi$  and $\Psi$ in $\F_G(A,B)$.\smallskip

\begin{remark}\label{ec-b-dist-r} For any $\Phi$  and $\Psi$ in $\F_G(A,B)$ the infimum in (\ref{ec-b-dist}) can be taken only over pure states $\omega$.
This follows from the freedom of choice of $R$, since the Bures distance between positive trace class operators defined in (\ref{B-d-s}) does not increase under partial trace: $\beta(\rho,\sigma)\geq \beta(\rho_A,\sigma_A)$ for any $\rho$ and $\sigma$ in $\T_+(\H_{AR})$ \cite{H-SCI,Watrous,Wilde}. \smallskip
\end{remark}

Inequality (\ref{B-d-s-r}) and definition (\ref{ecd}) imply that
\begin{equation}\label{D-B-in}
\frac{D_E^G(\Phi,\Psi)}{\sqrt{\|\Phi\|_{\diamond,E}^G}+\sqrt{\|\Psi\|_{\diamond,E}^G}}\leq \beta_E^G(\Phi,\Psi)\leq \sqrt{D_E^G(\Phi,\Psi)}
\end{equation}
for any maps $\Phi$  and $\Psi$ in $\F_G(A,B)$. \smallskip

\begin{lemma}\label{r-in+}
\emph{Let $\{\Phi_n\}$ be a sequence of maps in $\F_G(A,B)$. There exist a separable Hilbert space $\H_E$ and a sequence  $\{V_n\}$ of operators
in $\B_G(\H_A,\H_{BE})$ such that}
$$
\Phi_n(\rho)=\Tr_E V_n\rho V_n^*,\; \rho\in\T_G(\H_A),\quad  \textit{and}\quad \,\|V_{n+1}-V_{n}\|^G_E=\beta_E^G(\Phi_{n+1},\Phi_n)\;\;\;\forall n.
$$
\emph{If all the maps $\,\Phi_n$ lie in $\,\F_G^0(A,B)$ (correspondingly, in $\,\F(A,B)$) then all the operators $V_n\!$ can be taken in $\,\B^0_G(\H_A,\H_{BE})$ (correspondingly, in $\,\B(\H_A,\H_{BE})$).}
\end{lemma}\smallskip

\emph{Proof.} We may assume that all the maps $\Phi_n$ have representing operators $V_n$ in $\B_G(\H_A,\H_{BE})$, where $\H_E$
is an infinite-dimensional separable Hilbert space. Let $\H_{\tilde{E}}=\bigoplus_{k=1}^{+\infty} \H^k_E$, where $\H^k_E$ is a copy of $\H_E$ for each $k$. Let $\tilde{V}_{1}$ be an operator from $\D(\sqrt{G})$ into $\H_B\otimes\H_{\tilde{E}}=\bigoplus_{k=1}^{+\infty} \H_B\otimes\H^k_E$ defined by setting $\tilde{V}_{1}|\varphi\rangle=V_{1}|\varphi\rangle\oplus|0\rangle\oplus|0\rangle...$ for any $\varphi\in\D(\sqrt{G})$, where $V_{1}|\varphi\rangle\in \H_B\otimes\H^1_E$.

By Lemma \ref{r-in} below there is an operator $\tilde{V}_{2}$  from $\D(\sqrt{G})$ into $\H_B\otimes\H_{\tilde{E}}$ such that
$\|\tilde{V}_{2}-\tilde{V}_{1}\|^G_E=\beta_E^G(\Phi_2,\Phi_{1})$. This operator is  defined by setting
$$
\tilde{V}_{2}|\varphi\rangle=(I_B\otimes C_2)V_{2}|\varphi\rangle\oplus \left(I_B\otimes\sqrt{I_{E}-C_2^*C_2}\right)V_{2}|\varphi\rangle\oplus|0\rangle\oplus|0\rangle...
$$
for any $\varphi\in\D(\sqrt{G})$, where the first and the second summands here lie, respectively, in  $\H_B\otimes\H^1_E$ and $\H_B\otimes\H^2_E$ and $C_2$ is a contraction in $\B(\H_E)$.

Assume now that $\Phi=\Phi_2$, $\Psi=\Phi_3$ and $V_{\Phi}=(I_B\otimes C_2)V_{2}\oplus \left(I_B\otimes\sqrt{I_{E}-C_2^*C_2}\right)V_{2}$ is an operator
from $\D(\sqrt{G})$ into $\H_B\otimes (\H_E\oplus\H_E)$ representing the map $\Phi$. Then Lemma \ref{r-in}
below implies existence of an  operator $\tilde{V}_{3}$  from $\D(\sqrt{G})$ into $\H_B\otimes\H_{\tilde{E}}$ such that
$\|\tilde{V}_{3}-\tilde{V}_{2}\|^G_E=\beta_E^G(\Phi_3,\Phi_{2})$. This operator is  defined by setting
$$
\tilde{V}_{3}|\varphi\rangle=(I_B\otimes C_3 U_2)V_{3}|\varphi\rangle\oplus \left(I_B\otimes\sqrt{I_{E_2}-C_3^*C_3}\,U_2\right) V_{3}|\varphi\rangle\oplus|0\rangle\oplus|0\rangle...
$$
for any $\varphi\in\D(\sqrt{G})$, where the first and the second summands here lie, respectively, in  $\H_B\otimes(\H^1_E\oplus\H^2_E)$ and $\H_B\otimes(\H^3_E\oplus\H^4_E)$, $U_2$ is some unitary operator from $\H_E$ onto $\H_{E_2}\doteq\H_E\oplus\H_E$ and  $C_3$ is a contraction in $\B(\H_{E_2})$.

By using this way based on Lemma \ref{r-in} sequentially we obtain operators $\tilde{V}_{n}$, $n\geq 4$,  from $\D(\sqrt{G})$ into $\H_B\otimes\H_{\tilde{E}}$ such that
$\|\tilde{V}_{n+1}-\tilde{V}_{n}\|^G_E=\beta_E^G(\Phi_{n+1},\Phi_n)$ for all $n$. The operator $\tilde{V}_{n}$ is defined by setting
$$
\tilde{V}_{n}|\varphi\rangle=(I_B\otimes C_n U_{2^{n-2}})V_{n}|\varphi\rangle\oplus \left(I_B\otimes\sqrt{I_{E_{2^{n-2}}}\!-C_n^*C_n}\,U_{2^{n-2}}\right) V_{n}|\varphi\rangle\oplus|0\rangle\oplus|0\rangle...
$$
for any $\varphi\in\D(\sqrt{G})$, where the first and the second summands here lie, respectively, in  $\H_B\otimes(\H^1_E\oplus...\oplus\H^{2^{n-2}}_E)$ and $\H_B\otimes(\H^{2^{n-2}+1}_E\oplus...\oplus\H^{2^{n-1}}_E)$, $U_{2^{n-2}}$ is some unitary operator from $\H_E$ onto the direct sum $\H_{E_{2^{n-2}}}$ of $2^{n-2}$ copies of $\H_E$ and $C_n$ is a contraction in $\B(\H_{E_{2^{n-2}}})$.

By the construction $\Phi_n(\rho)=\Tr_{\tilde{E}} \tilde{V}_n\rho \tilde{V}_n^*$ for all $\rho\in\T_G(\H_A)$ and any $n$. \smallskip

It is easy to see that all the operators $\tilde{V}_{n}$ belong to the set  $\B_G(\H_A,\H_{B\tilde{E}})$. If all the operators $V_n\!$ lie in $\,\B^0_G(\H_A,\H_{BE})$ (correspondingly, in $\,\B(\H_A,\H_{BE})$) then all the operators $\tilde{V}_{n}$ lie in $\,\B^0_G(\H_A,\H_{B\tilde{E}})$ (correspondingly, in $\,\B(\H_A,\H_{B\tilde{E}})$). $\square$
% So, by redefining $\tilde{V}_{n}$ by $V_{n}$ and $\tilde{E}$ by $E$ we obtain the assertion of the lemma. $\square$

%\pagebreak

The construction from the proof of Theorem 1 in \cite{Kr&W} gives the following\smallskip

\begin{lemma}\label{r-in} \emph{Let $\Phi$ and $\Psi$ be any maps in $\F_G(A,B)$ with representing operators $V_{\Phi}$ and $V_{\Psi}$
in $\B_G(\H_A,\H_{BE})$, where $\H_E$ is a separable Hilbert space.  Let $\H_{\tilde{E}}=\H^1_E\oplus\H^2_E$, where $\H^1_E$ and $\H^2_E$ are copies of $\H_E$ and  $C$ is a contraction in  $\B(\H_{E}^2,\H_{E}^1)\cong\B(\H_{E})$.}\smallskip

\emph{The operators
$\tilde{V}_{\Phi}$ and $\tilde{V}_{\Psi}$  from $\D(\sqrt{G})\subseteq\H_A$ into $[\H_B\otimes\H_{E}^1]\oplus[\H_B\otimes\H_{E}^2]$ defined by setting
\begin{equation}\label{bar-v}
\!\tilde{V}_{\Phi}|\varphi\rangle=V_{\Phi}|\varphi\rangle\oplus|0\rangle,\quad
\tilde{V}^C_{\Psi}|\varphi\rangle=(I_B\otimes C)V_{\Psi}|\varphi\rangle\oplus \left(I_B\otimes\sqrt{I_{E}-C^*C}\right)V_{\Psi}|\varphi\rangle
\end{equation}
for any $\varphi\in\D(\sqrt{G})$ (where we assume that the operators $V_{\Phi}$ and $V_{\Psi}$ act from $\D(\sqrt{G})$ to $\H_B\otimes \H^1_E$ and
$\H_B\otimes \H^2_E$ correspondingly) are representing the maps $\Phi$ and $\Psi$.}\smallskip

\emph{There exists a contraction $C$ such that $\,\|\tilde{V}_{\Phi}-\tilde{V}^C_{\Psi}\|^G_E=\beta_E^G(\Phi,\Psi)\leq \sqrt{D_E^G(\Phi,\Psi)}$.}\smallskip
\end{lemma}

\emph{Proof.}  It is easy to see that the operators
$\tilde{V}_{\Phi}$ and $\tilde{V}_{\Psi}$ belong to the space $\B_G(\H_A,\H_{B\tilde{E}})$ and that they are representing the maps $\Phi$ and $\Psi$. By inequality (\ref{D-B-in}) to prove the lemma it suffices to show that $\inf_{C\in\B_1(\H_E)}\|\tilde{V}^C_{\Psi}-\tilde{V}_{\Phi}\|^G_E=\beta_E^G(\Phi,\Psi)$, where $\B_1(\H_E)$ is the unit ball of $\B(\H_E)$, and that this infimum is attainable.

We have
$$
\|\tilde{V}^C_{\Psi}-\tilde{V}_{\Phi}\|^G_E=\sup_{\rho\in\S_{G,E}}\sqrt{\Tr\shs\Phi(\rho)+\Tr\Psi(\rho)-2\Re\shs\Tr (I_B\otimes C)V_{\Phi}\rho V_{\Psi}^*},
$$
where $\S_{G,E}$ is the set of states $\rho$ in $\S(\H_A)$ such that $\Tr G\rho\leq E$ and  $V_{\Phi}\rho V_{\Psi}^*$ is the trace class operator well defined for any state $\rho$ in $\S_{G,E}$ by formula (\ref{ab-d}).

By Lemma \ref{qsl} the function $\rho\mapsto V_{\Phi}\rho V_{\Psi}^*$ is affine on the set $\S_{G,E}$ for any $E>0$ and takes values in $\T_+(\H_{BE})$.
This and the $\sigma$-weak compactness of the unit ball $\B_1(\H_E)$  (cf.\cite{B&R}) make it possible to apply
Ky Fan's minimax theorem (cf.\cite{Simons}) to change the
order of the optimization as follows\vspace{-5pt}
\begin{equation*}%\label{beta-d}
\begin{array}{rl}
\displaystyle\inf_{C\in\B_1(\H_E)}\|\tilde{V}^C_{\Psi}-\tilde{V}_{\Phi}\|^G_E &
=\displaystyle\inf_{C\in\B_1(\H_E)}\sup_{\rho\in\S_{{G},E}}\sqrt{\Tr\shs\Phi(\rho)+\Tr\Psi(\rho)-2\Re\shs\Tr(I_B\otimes C)V_{\Phi}\rho V_{\Psi} ^*}\\
&=\displaystyle\sup_{\rho\in\S_{{G},E}}\inf_{C\in\B_1(\H_E)} \sqrt{\Tr\shs\Phi(\rho)+\Tr\Psi(\rho)-2\Re\shs\Tr(I_B\otimes C)V_{\Phi}\rho V_{\Psi}^*}\\
&=\displaystyle\sup_{\rho\in\S_{{G},E}}\sqrt{\shs\Tr\shs\Phi(\rho)+\Tr\Psi(\rho)-2\sup_{C\in\B_1(\H_E)}|\Tr(I_B\otimes C)V_{\Phi}\rho V_{\Psi} ^*|}.
\end{array}
\end{equation*}
Note that all the above infima are attainable, since the expression under the first square root is a $\sigma$-weak continuous function of $C$ on the
$\sigma$-weak compact set $\B_1(\H_E)$.

By Lemma \ref{vbl} for any Hilbert space $\H_R$ the operators $V_{\Phi}\otimes I_R$ and $V_{\Psi}\otimes I_R$ are well defined on the set $\D(\sqrt{G}\otimes I_{R})$. Hence for any state $\rho$ in $\S_{{G},E}$ we have
 \begin{equation}\label{a-eq}
\!\sup_{C\in\B_1(\H_E)}|\Tr (I_B\otimes C)V_{\Phi}\rho V_{\Psi}^*|=\!\!\sup_{C\in\B_1(\H_E)}|\langle V_{\Psi}\otimes I_R\shs\varphi|I_{BR}\otimes C |V_{\Phi}\otimes I_R \shs\varphi\rangle|,\!\vspace{-5pt}
\end{equation}
where $\varphi$ is a purification of $\rho$, i.e. a vector in $\H_{A}\otimes\H_R$ such that $\Tr_R|\varphi\rangle\langle\varphi|=\rho$.

Since the vectors $V_{\Phi}\otimes I_R\shs|\varphi\rangle$ and $V_{\Psi}\otimes I_R\shs|\varphi\rangle$ in $\H_{BER}$ are purifications of the operators
$\Phi\otimes \id_R(\hat{\rho})$ and $\Psi\otimes \id_R(\hat{\rho})$ in $\T(\H_{BR})$, where $\hat{\rho}=|\varphi\rangle\langle\varphi|$, and the set $\B_1(\H_E)$ contains all isometries in $\B(\H_E)$, Uhlmann's theorem \cite{Uhlmann,Wilde} implies that the quantity in the r.h.s. of (\ref{a-eq}) is not less than the square root of the fidelity of these operators defined in (\ref{F-d-s}). On the other hand the r.h.s. of (\ref{a-eq}) is not greater than the square root of the fidelity, since
$\tilde{V}_{\Phi}\otimes I_R\shs|\varphi\rangle$ and $\tilde{V}^C_{\Psi}\otimes I_R\shs|\varphi\rangle$ are purifications of the operators
$\Phi\otimes \id_R(\hat{\rho})$ and $\Psi\otimes \id_R(\hat{\rho})$ as well and $\langle V_{\Psi}\otimes I_R\shs\varphi|I_{BR}\otimes C |V_{\Phi}\otimes I_R \shs\varphi\rangle=\langle \tilde{V}^C_{\Psi}\otimes I_R\shs\varphi|\tilde{V}_{\Phi}\otimes I_R \shs\varphi\rangle$.

Thus, since $\Tr\shs\Phi\otimes \id_R(\hat{\rho})=\Tr\shs\Phi(\rho)$ and  $\Tr\Psi\otimes \id_R(\hat{\rho})=\Tr\Psi(\rho)$, we have
$$
\begin{array}{rl}
\displaystyle\inf_{C\in\B_1(\H_E)}\|\tilde{V}^C_{\Psi}-\tilde{V}_{\Phi}\|^G_E &
=\displaystyle\sup_{\hat{\rho}_A\in\S_{{G},E}}\sqrt{\shs\Tr\shs\widehat{\Phi}(\hat{\rho})+\Tr\widehat{\Psi}(\hat{\rho})-
2\sqrt{F(\widehat{\Phi}(\hat{\rho}),\widehat{\Psi}(\hat{\rho}))}}\\
&=\displaystyle\sup_{\hat{\rho}_A\in\S_{{G},E}}\beta(\widehat{\Phi}(\hat{\rho}),\widehat{\Psi}(\hat{\rho}))=\beta_E^G(\Phi,\Psi),
\end{array}
$$
where $\widehat{\Theta}=\Theta\otimes \id_R$, $\Theta=\Phi,\Psi$,  the second equality follows from the definition (\ref{B-d-s}) of the Bures distance and the third one -- from Remark \ref{ec-b-dist-r}. $\square$
\medskip

\section{Properties of the cones $\F_G(A,B)$ and $\F^0_G(A,B)$}

\subsection{General properties and equivalent definition of $\,\F^0_G(A,B)$}

By definition a linear map $\Phi:\T_G(\H_A)\mapsto\T(\H_B)$ belongs to the cone $\F^0_G(A,B)$ if
there exist a separable Hilbert space $\H_E$ and an operator $V$ in $\B_G^0(\H_A,\H_{BE})$ such that
\begin{equation}\label{St-rep+}
\Phi(\rho)=\Tr_E V\rho V^*,\quad \rho\in\T_G(\H_A),
\end{equation}
where the operator $V\rho V^*$ is defined in (\ref{aa-d}). By Lemma \ref{qsl} the map $\Phi$ is uniformly continuous
on the set $\C_{G,E}=\{\rho\in\T_+(\H_A)\,|\,\Tr G\rho\leq E,\Tr\rho\leq 1\}$ for any $E>0$.\footnote{$\T_G(\H_A)$ is the linear span of all states $\rho\in\S(\H_A)$ with finite energy $\Tr G\rho$.}

A generalization of this property is presented in the following proposition, which also gives a characterization of the cone
$\F^0_G(A,B)$ in the case of discrete operator $G$. \smallskip

\begin{property}\label{ch-prop} \emph{Let $\H_R$ be a separable infinite-dimensional Hilbert space.}

A) \emph{If $\,\Phi$ is a map in $\F^0_G(A,B)$ then  the map $\Phi\otimes\id_R$ naturally defined on the set \footnote{$\,\T_G(\H_A)\otimes\T(\H_R)$ is the linear span of all operators $\rho\otimes\sigma$, where $\rho\in\T_G(\H_A)$ and $\sigma\in\T(\H_R)$.} $\,\T_G(\H_A)\otimes\T(\H_R)$ has a unique linear extension to the set $\T_{G\otimes I_R}(\H_{AR})$ (also denoted by $\Phi\otimes\id_R$) such that:}    \begin{enumerate}[a)]
         \item \emph{$\Phi\otimes\id_R(\omega)\geq 0$ for any positive operator $\omega$ in $\T_{G\otimes I_R}(\H_{AR})$;}
         \item \emph{the map $\Phi\otimes\id_R$ is uniformly continuous on the set $\,\C_{G\otimes I_R, E}\,$ for any $E>0$,  quantitatively,
        \begin{equation}\label{m-cb}
               \|\Phi\otimes\id_R(\omega_1)-\Phi\otimes\id_R(\omega_2)\|_1\leq 2\sqrt{\varepsilon \|\Phi\|^G_{\diamond,E}\|\Phi\|^G_{\diamond,4E/\varepsilon}}
         \end{equation}
             for any $\,\omega_1,\omega_2\in \C_{G\otimes I_R, E}$ such that $\|\omega_1-\omega_2\|_1\leq \varepsilon$.}\footnote{According to our notations $\,\C_{G\otimes I_R, E}\doteq\left\{\omega\in\T_+(\H_{AR})\,|\,\Tr G\omega_A\leq E,\,\Tr \omega\leq 1\right\}\,$. The r.h.s. of
             (\ref{m-cb}) tends to zero as $\varepsilon\rightarrow 0$ by condition (\ref{ph-cond}).}
        \end{enumerate}

B) \emph{If $G$ is a discrete operator (Def.1) and $\Phi:\T_G(\H_A)\mapsto\T(\H_B)$ is a linear map such that the map $\Phi\otimes\id_R$ has a  linear extension to the set $\,\T_{G\otimes I_R}(\H_{AR})$ possessing the above property a) and having continuous
restriction to the set $\,\C_{G\otimes I_R, E}$ for any $E>0$ then $\,\Phi$ belongs to the cone
$\F^0_G(A,B)$, i.e. it has representation (\ref{St-rep+}).}
\end{property}
\medskip

\emph{Proof.} A) Let $\Phi(\rho)=\Tr_E V\rho V^*$, where $V\in\B_{G}^0(\H_{A},\H_{BE})$.  In Section 3 (after Def.2) it is shown that for any separable Hilbert space $\H_R$ the map
$$
\T_{G\otimes I_R}(\H_{AR})\ni\omega\mapsto\Tr_E [V\otimes I_R]\omega[V\otimes I_R]^*
$$
is well defined and doesn't depend on $V$ (for given $\Phi$). Since $V\otimes I_R\in\B_{G\otimes I_R}^0(\H_{AR},\H_{BER})$ by Lemma \ref{vbl}, this
map naturally denoted by $\Phi\otimes \id_R$ belongs to the cone $\F^0_G(AR,BR)$. Since $\,\|V\otimes I_R\|^{G\otimes I_R}_{E}=\|V\|^G_{E}=\sqrt{\|\Phi\|^G_{\diamond,E}}$, Lemma \ref{qsl} implies continuity bound (\ref{m-cb}) showing 
uniform continuity of the map $\Phi\otimes \id_R$ on the set $\C_{G\otimes I_R,E}$ for any $E>0$.\smallskip

B) This assertion of the proposition follows from Theorem 3B in Section 5 (proved independently). $\square$ \smallskip

If  $V_1$ and $V_2$ are arbitrary operators in $\B_G(\H_A,\H_{BE})$, where $\H_E$ is any separable Hilbert space, then  Lemma \ref{qsl} shows that
the formula
\begin{equation}\label{St-c-r}
\Psi(\rho)= \Tr_E V_1\rho V_2^*
\end{equation}
where the operator $V_1\rho V_2^*$ is defined in (\ref{ab-d}), correctly defines a linear map form $\T_G(\H_A)$ into $\T(\H_B)$. If the operators  $V_1$ and $V_2$ lie in $\B^0_G(\H_A,\H_{BE})$ then Lemma \ref{qsl} implies that this map is uniformly continuous on the set $\C_{G,E}$ for any $E>0$. \smallskip

\begin{property}\label{l-span} \emph{A map $\,\Psi$ belongs to the linear span of $\,\F^0_G(A,B)$ (correspondingly, of $\,\F_G(A,B)$) if and only if it has representation (\ref{St-c-r}) for
some Hilbert space $\H_E$ and operators $V_1,V_2$ in $\B^0_G(\H_A,\H_{BE})$ (correspondingly, in $\B_G(\H_A,\H_{BE})$).}
\end{property}
\medskip

\emph{Proof.} If $\Psi$ belongs to the linear span of $\,\F^0_G(A,B)$ then it can be represented as follows
$$
\Psi=(\Phi_1-\Phi_2)+i(\Phi_3-\Phi_4),\quad \Phi_k\in\F^0_G(A,B),\;\; k=1,2,3,4.
$$
Let $V_k\in\B^0_G(\H_A,\H_{BE_k})$ be a representing operator for the map $\Phi_k$. Then the operators
$W_1$ and $W_2$ from $\D(\sqrt{G})$ to $\H_{B\tilde{E}}\doteq \H_{B}\otimes(\H_{E_1}\oplus\H_{E_2}\oplus\H_{E_3}\oplus\H_{E_4})$ defined by settings
$$
W_1 |\varphi\rangle=V_1 |\varphi\rangle\oplus -V_2 |\varphi\rangle\oplus i V_3 |\varphi\rangle\oplus -i V_4 |\varphi\rangle,\;\;
W_2 |\varphi\rangle=V_1 |\varphi\rangle\oplus V_2 |\varphi\rangle\oplus V_3 |\varphi\rangle\oplus V_4 |\varphi\rangle
$$
for any $\varphi\in\D(\sqrt{G})$, belong to the space $\B^0_G(\H_A,\H_{B\tilde{E}})$. It is easy to see that $\Psi(\rho)= \Tr_{\tilde{E}} W_1\rho W_2^*$
for any $\rho\in\T_G(\H_A)$.\smallskip

If $\Psi$ has representation (\ref{St-c-r}) then $\Psi=\frac{1}{4}(\Phi_1-\Phi_2+i\Phi_3-i\Phi_4)$, where
$$
\Phi_1(\rho)= \Tr_{E} (V_1+V_2)\shs\rho\shs (V_1+V_2)^*,\;\quad\Phi_2(\rho)= \Tr_{E} (V_1-V_2)\shs\rho\shs (V_1-V_2)^*,
$$
$$
\Phi_3(\rho)= \Tr_{E} (V_1+iV_2)\shs\rho\shs (V_1+iV_2)^*,\quad\Phi_4(\rho)= \Tr_{E} (V_1-iV_2)\shs\rho\shs (V_1-iV_2)^*
$$
are maps in $\,\F^0_G(A,B)$. \smallskip

The case of the cone $\F_G(A,B)$ is considered similarly. $\square$

\subsection{Kraus representation}

A CP linear map $\Phi:\T(\H_A)\rightarrow\T(\H_B)$ is characterized by the Kraus representation
\begin{equation}\label{K-rep}
\Phi(\rho)=\sum_k V_k\rho V^*_k, \quad \rho\in \T(\H_A),
\end{equation}
where $\{V_k\}$ is a finite or countable collection of operators in $\B(\H_A,\H_B)$ such that $\|\sum_k  V^*_kV_k\|=\|\Phi\|_{\diamond}$ \cite{H-SCI,Watrous,Wilde}. A similar characterizations
of the cone $\F_G(A,B)$ and its subcone
$\F^0_G(A,B)$ are presented in the following \smallskip

\begin{property}\label{Kraus} A) \emph{A map $\Phi$ belongs to the cone $\F_G(A,B)$ if and only if it can be represented as
\begin{equation}\label{K-rep}
\Phi(\rho)=\sum_k V_k\rho V^*_k, \quad \rho\in \T_G(\H_A),
\end{equation}
where $\{V_k\}$ is a finite or countable collection of operators in $\B_G(\H_A,\H_B)$ such that \footnote{If condition (\ref{K-c}) holds for some $E>0$ then it holds for all $E>0$. This follows from concavity of the function $E\mapsto\|\{V_k\}\|_E^G$ on $\mathbb{R}_+$.}
\begin{equation}\label{K-c}
\|\{V_k\}\|_E^G\doteq \sup_{\rho\in\S(\H_A):\Tr G\rho\leq E}\sum_k \Tr V_k\rho V^*_k<+\infty\;\;\textit{ for some  }\;\;E>0.
\end{equation}
The quantity $\|\{V_k\}\|_E^G$ coincides with the ECD norm $\|\Phi\|_{\diamond,E}^G$ defined in (\ref{E-sn})}.\footnote{$\T_G(\H_A)$ is the linear span of all the states $\rho$ with finite energy $\Tr G\rho$.} \smallskip

B) \emph{A map $\Phi$ belongs to the cone $\F^0_G(A,B)$ if and only if it has representation (\ref{K-rep}), where $\{V_k\}$ is collection of operators in $\B^0_G(\H_A,\H_B)$  such that}
\begin{equation}\label{K-c+}
\|\{V_k\}\|_E^G=o\shs(E)\;\;\textit{ as  }\;\;E\rightarrow+\infty.
\end{equation}
\end{property}
\medskip
\begin{remark}\label{Kraus-r}
Property (\ref{K-c+}) implies that $V_k\in\B^0_G(\H_A,\H_B)$ for all $k$. This follows from the condition (\ref{s-cond}) characterizing
$\sqrt{G}$-infinitesimal operators, since it is easy to see that $\|V_k\|_E^G\leq\sqrt{ \|\{V_k\}\|_E^G}\,$ for each $k$. This condition also shows that
property (\ref{K-c+}) holds for any \emph{finite} collection $\{V_k\}\subset\B^0_G(\H_A,\H_B)$.
\medskip
\end{remark}

\emph{Proof.} If $\Phi\in\F_G(A,B)$ then there exist a Hilbert space $\H_E$ and an operator $V$ in $\B_G(\H_A,\H_{BE})$ such that
$\Phi(\rho)=\Tr_E V\rho V^*$ for all $\rho$ in $\T_G(\H_A)$.

Let $\{\tau_k\}$ be an orthonormal basis in $\H_E$ and $\H_R\cong\H_A$. Lemma \ref{vbl} implies that the operator $V\otimes I_R$ belongs to the space $\B_G(\H_{AR},\H_{BER})$. By Lemma \ref{SK-l} below for each $k$ there exist unique operators $V_k\in\B_G(\H_A,\H_B)$ and $\widehat{V}_k\in\B_{G\otimes I_R}(\H_{AR},\H_{BR})$ such that
$$
\langle\psi\otimes\tau_k|V|\varphi\rangle=\langle\psi|V_{k}|\varphi\rangle\quad\textup{and}\quad \langle\theta\otimes\tau_k|V\otimes I_R|\eta\rangle=\langle\theta|\widehat{V}_{k}|\eta\rangle
$$
for all $\varphi\in\D(\sqrt{G})$, $\psi\in\H_B$, $\eta\in\D(\sqrt{G}\otimes I_R)$ and $\theta\in\H_{BR}$.
These relations imply that $\widehat{V}_{k}=V_k\otimes I_R$ and hence
\begin{equation}\label{Tr-eq}
 \! \Tr [V\otimes I_R]|\eta\rangle\langle\eta| [V\otimes I_R]^*=\sum_k \Tr [V_k\otimes I_R]|\eta\rangle\langle\eta|[V_k\otimes I_R]^* \quad\forall\eta\in\D(\sqrt{G}\otimes I_R).
\end{equation}
It follows that $[V\otimes I_R]|\eta\rangle=\sum_k [V_k\otimes I_R]|\eta\rangle\otimes |\tau_k\rangle$ for any $\eta\in\D(\sqrt{G}\otimes I_R)$ (the convergence of the last series follows from the equality
$\sum_k\|[V_k\otimes I_R]\eta\|^2=\|[V\otimes I_R]\eta\|^2$  which is a partial case of (\ref{Tr-eq})). Thus, for any $\eta\in\D(\sqrt{G}\otimes I_R)$ we have
\begin{equation}\label{br}
[V\otimes I_R]|\eta\rangle\langle\eta| [V\otimes I_R]^*=\sum_{k,j} [V_k\otimes I_R]|\eta\rangle\langle\eta| [V_j\otimes I_R]^*\otimes |\tau_k\rangle\langle\tau_j|.
\end{equation}
Since any state $\rho$ in $\S_G(\H_A)$ can be represented as $\Tr_R|\eta\rangle\langle\eta|$ for some unit vector $\eta$ in $\D(\sqrt{G}\otimes I_R)$,
by taking partial trace over the spaces $\H_E$ and $\H_R$ in the above equality we obtain (\ref{K-rep}). It follows that the r.h.s. of (\ref{K-c})
coincides with $[\|V\|_E^G]^2=\|\Phi\|_{\diamond,E}^G$.

Let $\Phi$ be a map with representation (\ref{K-rep}), $\H_R\cong\H_A$ and $\H_E$ a Hilbert space whose dimension coincides with the cardinality
of the set $\{V_k\}$. Lemma \ref{vbl} implies that for each $k$ the operator $V_k\otimes I_R$ belongs to the space $\B_G(\H_{AR},\H_{BR})$.
Since condition (\ref{K-c}) implies that
$\sum_k\|[V_k\otimes I_R] \eta\|^2<+\infty$  for any $\eta\in\D(\sqrt{G}\otimes I_R)$,  we may define for given basis $\{\tau_k\}$ in $\H_E$ the operator
$\widehat{V}$ from $\D(\sqrt{G}\otimes I_R)$ into $\H_{BER}$ by setting $\widehat{V}|\eta\rangle=\sum_k [V_k\otimes I_R]|\eta\rangle\otimes |\tau_k\rangle$ for any $\eta\in\D(\sqrt{G}\otimes I_R)$. It is easy to see that $\widehat{V}=V\otimes I_R$, where $V$ is the operator from $\D(\sqrt{G})$ into $\H_{BE}$ defined by setting $V|\varphi\rangle=\sum_k V_k|\varphi\rangle\otimes |\tau_k\rangle$ for any $\varphi\in\D(\sqrt{G})$.
This  implies that equality (\ref{br}) holds for any $\eta\in\D(\sqrt{G}\otimes I_R)$. Since any state $\rho$ in $\S_G(\H_A)$ can be represented as $\Tr_R|\eta\rangle\langle\eta|$ for some unit vector $\eta\in\D(\sqrt{G}\otimes I_R)$,
by taking partial trace over the spaces $\H_E$ and $\H_R$ in this equality we obtain $\Phi(\rho)=\Tr_E V\rho V^*$ for all $\rho$ in $\S_G(\H_A)$.\smallskip

B) This assertion follows from the proof of part A, where it is shown that the r.h.s. of (\ref{K-c})
coincides with $[\|V\|_E^G]^2=\|\Phi\|_{\diamond,E}^G$. It suffices only to note that
$\|V_k\|_E^G\leq\|\{V_k\}\|_E^G$ for each $k$ and to use the conditions (\ref{s-cond}) and (\ref{ph-cond}) characterizing
$\sqrt{G}$-infinitesimal operators and the maps in $\F^0_G(A,B)$ correspondingly. $\square$ \medskip

\begin{lemma}\label{SK-l}\emph{If $\,V$ is an operator in $\B_G(\H_A,\H_{BE})$ then
for any unit vector $\tau$ in $\H_E$ there is a unique operator $V_{\tau}$ in $\B_G(\H_A,\H_{B})$ such that $\|V_{\tau}\|_E^G\leq \|V\|_E^G$ for any $E>0$ and}
\begin{equation}\label{SK-1-r}
\langle\psi\otimes\tau|V|\varphi\rangle=\langle\psi|V_{\tau}|\varphi\rangle\quad\textit{ for all }\varphi\in\D(\sqrt{G})\textit{ and }\psi\in\H_B.
\end{equation}
\end{lemma}

\emph{Proof.} Denote by $\hat{\H}_A$ the Hilbert space obtained by equipping the linear space $\D(\sqrt{G})$
with the inner product (\ref{in-p}) for some $E>0$. Then the l.h.s. of (\ref{SK-1-r}) is a sesquilinear form on $\H_B\times\hat{\H}_A$
bounded above by $\sn V\sn^{G}_{E}\|\psi\|\|\varphi\|^{G}_{E}$, where $\sn V\sn^{G}_{E}$ is the norm of $V$ as an operator from $\hat{\H}_A$ to $\H_{BE}$ defined in (\ref{eq-norms-2}) and $\|\varphi\|^{G}_{E}=\sqrt{\|\varphi\|^2+\|\sqrt{G}\varphi\|^2/E}$ is the norm of $\varphi\in\hat{\H}_A$. So, there is a unique bounded linear operator $V_{\tau}$ from $\hat{\H}_A$ into $\H_B$ with $\sn V_{\tau}\sn^{G}_{E}\leq\sn V\sn^{G}_{E}$ such that this sesquilinear form is equal to $\langle\psi|V_{\tau}|\varphi\rangle$ for any $\varphi\in\hat{\H}_A$ and $\psi\in\H_B$. Since $\hat{\H}_A$ and $\D(\sqrt{G})$ coincide as linear spaces,  the operator $V_{\tau}$ can be considered as a $\sqrt{G}$-bounded operator from $\H_A$ into $\H_B$ which does not depend on the chosen value of $E>0$. It follows that $\sn V_{\tau}\sn^{G}_{E}\leq\sn V\sn^{G}_{E}$ for any $E>0$. Hence, the second relation in (\ref{s-tr}) implies that
$\|V_{\tau}\|^{G}_{E}\leq\|V\|^{G}_{E}$ for any $E>0$. $\square$

\subsection{Generalized version of the Kretschmann-Schlingemann-Werner theorem}

The Kretschmann-Schlingemann-Werner theorem  (the KSW-theorem in what follows) obtained in \cite{Kr&W} states that
\begin{equation*}
 \frac{\|\Phi-\Psi\|_{\diamond}}{\sqrt{\|\Phi\|_{\diamond}}+\sqrt{\|\Psi\|_{\diamond}}}\leq\inf_{V_{\Phi},V_{\Psi}}\|V_{\Phi}-V_{\Psi}\|=\beta(\Phi,\Psi) \leq \sqrt{\|\Phi-\Psi\|_{\diamond}},
\end{equation*}
for any maps $\Phi$ and $\Psi$ in the cone $\F(A,B)$, where the infimum is over all common Stinespring representations
\begin{equation}\label{c-S-r}
\Phi(\rho)=\Tr_E V_{\Phi}\rho V^*_{\Phi},\qquad\Psi(\rho)=\Tr_E V_{\Psi}\rho V^*_{\Psi}
\end{equation}
and $\beta(\Phi,\Psi)$ is the Bures distance between the maps $\Phi$ and $\Psi$ defined in (\ref{b-dist}). \smallskip

Lemmas \ref{l-in} and \ref{r-in} imply the following generalized version of the KSW theorem.\smallskip

\begin{theorem}\label{KSW++} \emph{Let $G$ be a positive operator on $\H_A$ satisfying condition (\ref{G-cond}), $\|\!\cdot\!\|^G_E$, $\|\!\cdot\!\|^G_{\diamond,E}$ and $\beta_E^G$ the operator E-norm, the ECD norm and the energy constrained  Bures distance defined, respectively, in (\ref{ec-on}), (\ref{E-sn}) and (\ref{ec-b-dist}). Let $\H_E$ be an infinite-dimensional separable Hilbert space.}
\emph{For any maps  $\Phi$ and $\Psi$ in $\F_G(A,B)$ the following relations hold
$$
\frac{\|\Phi-\Psi\|_{\diamond,E}^G}{\sqrt{\|\Phi\|_{\diamond,E}^G}+\sqrt{\|\Psi\|_{\diamond,E}^G}}\leq\inf_{V_{\Phi},V_{\Psi}}\|V_{\Phi}-V_{\Psi}\|^{G}_E\leq \beta_E^G(\Phi,\Psi)\leq \sqrt{\|\Phi-\Psi\|_{\diamond,E}^G},
$$
where the infimum is over all operators $V_{\Phi}$ and $V_{\Psi}$ in $\B_G(\H_A,\H_{BE})$ representing the maps $\Phi$ and $\Psi$, i.e. such that (\ref{c-S-r}) holds for any state $\rho$ in $\,\S(\H_A)$ with finite $\Tr G\rho$.}
\end{theorem}\smallskip

\textbf{Note:} In contrast to the original KSW theorem mentioned above  and to the\break  \emph{E}-version of this theorem presented in \cite[Section 3]{ECN}, Theorem \ref{KSW++} \emph{do not assert} that $\inf_{V_{\Phi},V_{\Psi}}\|V_{\Phi}-V_{\Psi}\|^{G}_E=\beta_E^G(\Phi,\Psi)$ for any maps  $\Phi$ and $\Psi$ in $\F_G(A,B)$.

\section{On completion of the set of Hermitian-preserving completely bounded linear maps w.r.t. the ECD norm}
Let  $\Y(A,B)$ be the real linear space of all Hermitian-preserving completely bounded linear maps from $\T(\H_A)$ into $\T(\H_B)$.
The space $\Y(A,B)$ endowed with the diamond norm (\ref{d-n-def}) is a real Banach space, but $\Y(A,B)$ is not complete w.r.t. the ECD norm (\ref{E-sn}) if $G$ is an unbounded operator. In this section we describe the real Banach space $\Y_G(A,B)$ containing the completion of $\Y(A,B)$  w.r.t. the ECD norm, which coincides with this completion if $G$ is a discrete unbounded operator (Def.1).\smallskip

Let $\Y_G(A,B)$ be the set of all linear maps $\Phi$ from the subset $\T_G(\H_A)\subset\T(\H_A)$ into $\T(\H_B)$  with the following properties\footnote{$\T_G(\H_A)$ is the linear span of all states $\rho$ in $\S(\H_A)$ with finite energy $\Tr G\rho$.}
\begin{enumerate}[1)]
  \item $\Phi(\rho^*)=[\Phi(\rho)]^*$ for all $\rho$ in $\T_G(\H_A)$, i.e. $\Phi$ is Hermitian preserving;
  \item For any separable Hilbert space $\H_R$ the map $\Phi\otimes\id_R$ naturally defined on the set $\,\T_G(\H_A)\otimes\T(\H_R)$ has a linear extension to the set $\T_{G\otimes I_R}(\H_{AR})$, which is continuous on the set $\,\C_{G\otimes I_R, E}\doteq\left\{\omega\in\T_+(\H_{AR})\,|\,\Tr G\omega_A\leq E,\,\Tr \omega\leq 1\right\}\,$ for any $E>0$.
\end{enumerate}

Show first that the extension of $\Phi\otimes\id_R$ mentioned in property 2 is unique. Assume that $\Theta$ and $\Theta'$
are extensions of $\Phi\otimes\id_R$, which are continuous on the set $\,\C_{G\otimes I_R, E}$ for any $E>0$. It suffices to show that
$\Theta$ and $\Theta'$ coincide at a pure state $|\eta\rangle\langle\eta|$, where $\eta$ is a any unit vector in $\D(\sqrt{G}\otimes I_R)$ (since any state in $\S_{G\otimes I_R}(\H_{AR})$ is decomposed into a convex combination of such pure states).

Any unit vector $\eta$  in $\D(\sqrt{G}\otimes I_R)$ has the representation
$$
|\eta\rangle=\sum_{i}|\varphi_i\rangle\otimes|\psi_i\rangle,
$$
where $\{\varphi_i\}$ and $\{\psi_i\}$ are collections of vectors in $\D(\sqrt{G})$ and $\H_R$ such that  $\sum_{i}\|\sqrt{G}\varphi_i\|^2=E<+\infty$ and $\langle\psi_i|\psi_j\rangle=\delta_{ij}$.
For any given $n$ let
$$
|\eta_n\rangle\langle\eta_n|=\sum_{i,j=1}^n|\varphi_i\rangle\langle\varphi_j|\otimes|\psi_i\rangle\langle\psi_j|.
$$
Then
$$
\Theta(|\eta_n\rangle\langle\eta_n|)=
\sum_{i,j=1}^n\Phi(|\varphi_i\rangle\langle\varphi_j|)\otimes|\psi_i\rangle\langle\psi_j|=\Theta'(|\eta_n\rangle\langle\eta_n|).
$$
Since the state $|\eta\rangle\langle\eta|$ and all the operators $|\eta_n\rangle\langle\eta_n|$ belong to the set $\,\C_{G\otimes I_R, E}$, the continuity of the maps $\Theta$ and $\Theta'$ on this set  implies that $\Theta(|\eta\rangle\langle\eta|)=
\Theta'(|\eta\rangle\langle\eta|)$.\smallskip

Let $\Y^+_G(A,B)$ be the subset of $\Y_G(A,B)$ consisting of maps $\Phi$ such that $\Phi\otimes\id_R(\omega)$ is positive
for any positive $\omega\in\T_{G\otimes I_R}(\H_{AR})$, where $R$ is any system. It is clear that $\Y^+_G(A,B)$ is a cone in $\Y_G(A,B)$.
%The observations in Section 4.1 show that $\Y^+_G(A,B)$ contains the cone $\F_G(A,B)$.

It is easy to see that $\Y_G(A,B)$ is a real linear space and that $\|\cdot\|_{\diamond,E}^G$ is a norm on  $\Y_G(A,B)$ for any $E>0$.
By repeating the arguments in \cite{W-EBN} one can show  that for any given $\Phi\in\Y_G(A,B)$ the nondecreasing nonnegative function $E\mapsto\|\Phi\|_{\diamond,E}^G$ is concave on $\mathbb{R}_+$. It implies relations (\ref{ECD-n-eq}) which show the equivalence of all the
norms $\|\cdot\|_{\diamond,E}^G$, $E>0$, on the set $\Y_G(A,B)$ for given operator $G$.  In what follows we will assume that the space $\Y_G(A,B)$ is endowed with the norm $\|\cdot\|_{\diamond,E}^G$ for some $E>0$.\smallskip

\begin{theorem}\label{third}
A) \emph{$\Y_G(A,B)$ is a real Banach space containing the set $\,\Y(A,B)$.}\smallskip

B) \emph{If $G$ is a discrete operator (Def.1) then $\,\Y_G(A,B)$ is the completion of $\,\Y(A,B)$ w.r.t. the norm $\|\!\cdot\!\|_{\diamond,E}^G$ and
$\,\Y^+_G(A,B)=\F^0_G(A,B)$.}\footnote{The cone $\F^0_G(A,B)$ is defined in Section 3, where it is shown that it coincides with the completion of the cone $\,\F(A,B)$ of CP maps in $\,\Y(A,B)$ w.r.t. the norm $\|\!\cdot\!\|_{\diamond,E}^G$.}
\end{theorem}\smallskip

\emph{Proof.} A) Note first  that condition (\ref{G-cond}) implies that
\begin{equation}\label{ECDN-d+}
  \|\Phi\otimes \id_R(\omega)\|_1\leq \|\Phi\|^G_{\diamond,E}
\end{equation}
for any $\Phi$ in $\Y_G(A,B)$ and any operator $\omega$ in $\T_+(\H_{AR})$ such that $\Tr G\omega_A\leq E$ and $\Tr \omega\leq 1$, where $\H_R$ is a separable Hilbert space. Indeed, let $\omega$ be such an operator and $r=\Tr\omega$. Then $\hat{\omega}\doteq r^{-1}\omega$ is a state such that
$\,\Tr {G}\hat{\omega}_A\leq E/r$. So, by using concavity of the function $E\mapsto\|\Phi\|^G_{\diamond,E}$ on $\mathbb{R}_+$ mentioned above and Lemma \ref{WL} we obtain
$$
\|\Phi\otimes \id_R(\omega)\|_1=r\|\Phi\otimes \id_R(\hat{\omega})\|_1\leq r\|\Phi\|^G_{\diamond,E/r}\leq \|\Phi\|^G_{\diamond,E}.
$$

Let $\{\Phi_n\}$ be a  Cauchy  sequence in $\,\Y_G(A,B)$ and $R$ an infinite-dimensional quantum system. Then for any operator $\omega$ in $\T_{G\otimes I_R}(\H_{AR})$ the sequences $\{\Phi_n\otimes\id_R(\omega)\}$ and $\{\Phi_n(\omega_A)\}$  are  Cauchy  sequences in $\T(\H_{BR})$ and $\T(\H_{B})$ correspondingly. Hence they have limits which will be denoted, respectively,  by $\Theta(\omega)$ and $\Phi(\omega_A)$.
By this way we define the Hermitian-preserving linear maps $\Theta:\T_{G\otimes I_R}(\H_{AR})\rightarrow\T(\H_{BR})$ and $\Phi:\T_{G}(\H_{A})\rightarrow\T(\H_{B})$.

Note that
\begin{equation}\label{lr-1}
\lim_{n\rightarrow+\infty}\,\sup_{\omega\in\C_{G\otimes I_R ,E}}\|\Phi_n\otimes\id_R(\omega)-\Theta(\omega)\|_1=0
\end{equation}
(where $\,\C_{G\otimes I_R, E}\doteq\left\{\omega\in\T_+(\H_{AR})\,|\,\Tr G\omega_A\leq E,\,\Tr \omega\leq 1\right\}$). Indeed, if this relation does not
hold then (by passing to a subsequence)
we may assume that there is $\varepsilon>0$ and a sequence $\{\omega_n\}\subset\C_{G\otimes I_R ,E}$ such that  $\|\Phi_n\otimes\id_R(\omega_n)-\Theta(\omega_n)\|_1\geq\varepsilon$. By choosing $n$ such that $\|\Phi_n-\Phi_m\|_{\diamond,E}^G<\varepsilon/2$ for all $m>n$, we have
$$
\|\Phi_n\otimes\id_R(\omega_n)-\Theta(\omega_n)\|_1\leq \|\Phi_n\otimes\id_R(\omega_n)-\Phi_m\otimes\id_R(\omega_n)\|_1+\|\Phi_m\otimes\id_R(\omega_n)-\Theta(\omega_n)\|_1.
$$
It follows from (\ref{ECDN-d+}) that the first term in the r.h.s. of this inequality does not exceed $\|\Phi_n-\Phi_m\|_{\diamond,E}^G<\varepsilon/2$, while the second one can be made less than
$\varepsilon/2$ by choosing sufficiently large $m$. This contradicts the above assumption.

By definition of  $\Y_G(A,B)$ all the functions $\,\omega\mapsto\Phi_n\otimes\id_R(\omega)\,$ are continuous on the set $\C_{G\otimes I_R, E}$ for each $E>0$. So, relation (\ref{lr-1}) implies that the function $\omega\mapsto\Theta(\omega)$ is also continuous on the set $\C_{G\otimes I_R, E}$ for each $E>0$.

By the definitions of $\Theta$ and $\Phi$ the map $\Theta$ coincides with the map $\Phi\otimes\id_R$ on the set $\,\T_G(\H_A)\otimes\T(\H_B)$.
Thus, $\Theta$ is the extension of $\Phi\otimes\id_R$ mentioned in property 2 of the above definition of $\Y_G(A,B)$.

Relation (\ref{lr-1}) implies that $\|\Phi_n-\Phi\|_{\diamond,E}^G$ tends to zero as $n\rightarrow+\infty$.\smallskip

B) We have to show that any map $\Phi$ in  $\Y_G(A,B)$ can be approximated by a sequence of maps in  $\Y(A,B)$ w.r.t. the norm $\|\!\cdot\!\|_{\diamond,E}^G$ provided that the operator $G$ has representation (\ref{H-rep}).

For given natural $n$ denote by $\H_A^n$ the linear span of the vectors $\tau_0,...,\tau_{n-1}$, i.e. $\H_A^n$ is the subspace of $\H_A$ corresponding to the minimal $n$  eigenvalues of ${G}$ (taking the multiplicity into account). Denote by $P_n$ the projector onto $\H_A^n$. Consider the quantum channel $\Pi_n(\rho)=P_n\rho P_n+[\Tr \bar{P}_n\rho]|\tau_0\rangle\langle\tau_0|$, where $\bar{P}_n=I_{A}-P_n$.

Let $\Phi$ be a map in  $\Y_G(A,B)$. Note first that  the map $\Phi_n=\Phi\circ\Pi_n$ belongs to the set $\Y(A,B)$ for  any natural $n$. Indeed, for arbitrary state $\omega$ in $\S(\H_{AR})$ we have
$$
\|\Phi_n\otimes\id_R(\omega)\|_1=\|\Phi\otimes\id_R(\Pi_n\otimes\id_R(\omega))\|_1\leq \|\Phi\|_{\diamond,E_n}^G,
$$
where the inequality follows from definition of the norm $\|\!\cdot\!\|_{\diamond,E}^G$, since $\Tr G \Pi_n(\rho)\leq E_n$ for any $\rho\in\S(\H_A)$. So, it follows from  the definition of the diamond norm that
$$
\|\Phi_n\|_{\diamond}\leq \|\Phi\|_{\diamond,E_n}^G<+\infty\quad \forall n.
$$

Let $\tilde{G}$ be a positive operator on $\H_R$ isomorphic to the operator $G$. By Lemma \ref{ecdn-sd} below we have
\begin{equation}\label{dn-e}
\|\Phi-\Phi_n\|_{\diamond,E}^G=\sup_{\omega\in\S_{G,\tilde{G},E}}\|\Phi\otimes\id_R(\omega-\Pi_n\otimes\id_R(\omega))\|_1,
\end{equation}
where $\S_{G,\tilde{G},E}\doteq\{\omega\in\S(\H_{AR})\,|\,\Tr G\omega_A\leq E, \Tr \tilde{G}\omega_R\leq E\}$. \smallskip

Let $\bar{P}_n=I_A-P_n$.  By using the inequality
$$
\|\bar{P}_n\otimes I_R\cdot\omega\cdot P_n\otimes I_R\|_1\leq \sqrt{\Tr\bar{P}_n\otimes I_R\;\omega}
$$
easily proved for any $\omega\in\S(\H_{AR})$ via the operator Cauchy-Schwarz inequality (see the proof of Lemma 11.1 in \cite{H-SCI}) and by noting that $\,\Tr G\rho\leq E\,$ implies $\,\Tr \bar{P}_n\rho\leq E/E_n\,$  for any $\rho\in\S(\H_A)$ we obtain
\begin{equation}\label{norm-est}
\begin{array}{c}
\|\omega-\Pi_n\otimes\id_R(\omega)\|_1\leq 2\|P_n\otimes I_R\cdot\omega\cdot \bar{P}_n\otimes I_R\|_1+\|\bar{P}_n\otimes I_R\cdot\omega\cdot\bar{P}_n\otimes I_R\|_1\\\\
+\|\Tr_A\bar{P}_n\otimes I_R\shs\omega\|_1\leq 2\sqrt{\Tr\bar{P}_n\omega_A}+2\Tr\bar{P}_n\omega_A \leq4\sqrt{\Tr\bar{P}_n\omega_A}\leq4\sqrt{E/E_n}
\end{array}
\end{equation}
for any state $\omega$ in $\S(\H_{AR})$ such that $\Tr G\omega_A\leq E$.

By the Lemma in \cite{H-c-w-c} the set of states $\rho$ satisfying the inequality $\Tr G\rho\leq E$ is compact for any $E>0$. So, Corollary 6 in \cite{AQC} implies that the set $\S_{G,\tilde{G},E}$ in (\ref{dn-e}) is a compact subset of $\C_{G\otimes I_R, E}$  for any $E>0$. It follows that the continuous function $\omega\mapsto\Phi\otimes\id_R(\omega)$
is uniformly continuous on $\S_{G,\tilde{G},E}$. Thus, estimate (\ref{norm-est}) implies that the r.h.s. of (\ref{dn-e}) tends to zero as $n\rightarrow+\infty$.

To prove that $\,\Y^+_G(A,B)=\F^0_G(A,B)$ we have to show that any map $\Phi$ in $\Y^+_G(A,B)$ lies in $\F^0_G(A,B)$.
We will use the sequence $\{\Phi_n\}$ constructed before. In this case it consists of  maps in $\F(A,B)$. Since the cone $\F^0_G(A,B)$ is complete
w.r.t. the ECD norm by Theorem \ref{main}, the limit map $\Phi$ of the sequence $\{\Phi_n\}$ belongs to this cone. $\square$
\smallskip

\begin{lemma}\label{ecdn-sd} \emph{Let $\tilde{G}$ be a positive operator on $\H_R$ isomorphic to the operator $G$. Then
$$
\|\Phi\|_{\diamond,E}^G=\sup_{\omega\in\S_{G,\tilde{G},E}}\|\Phi\otimes\id_R(\omega)\|_1
$$
for any $\Phi$ in $\Y_G(A,B)$, where $\S_{G,\tilde{G},E}\doteq\{\omega\in\S(\H_{AR})\,|\,\Tr G\omega_A\leq E, \Tr \tilde{G}\omega_R\leq E\}$.}
\end{lemma}\medskip

\emph{Proof.} Since the system $R$ in definition (\ref{E-sn}) is assumed arbitrary, the supremum in (\ref{E-sn}) can be taken over all pure states $\omega$ in $\S(\H_{AR})$ satisfying the condition $\,\Tr {G}\omega_A\leq E$. Since for any such pure state $\omega$ the partial states $\omega_A$ and $\omega_R$  have the same nonzero spectrum,  by applying local partial isometry transformation of the system $R$ this state can be transformed into a state  $\omega'$ belonging to the set $\S_{G,\tilde{G},E}$. It suffices to note that $\,\|\Phi\otimes\id_R(\omega)\|_1=\|\Phi\otimes\id_R(\omega')\|_1$.  $\square$ \smallskip

It is well known that any map $\Phi$ in $\Y(A,B)$ can be represented as $\Psi_1-\Psi_2$, where $\Psi_1$ and $\Psi_2$ are maps in $\F(A,B)$ \cite{Paul,Watrous}. Theorem \ref{third}B gives a reason for the following\smallskip

\textbf{Conjecture.} Any map $\Phi$ in $\Y_G(A,B)$, where $G$ is a positive discrete unbounded operator, can be represented as $\Phi=\Psi_1-\Psi_2$, where $\Psi_1$ and $\Psi_2$  are maps in $\F^0_G(A,B)$.\footnote{I would be grateful for any comments concerning this conjecture.}

\smallskip

Theorem \ref{third}B implies the following characterisation of the class of maps in $\Y_G(A,B)$ which have such representation.
\smallskip

\begin{property}\label{p-c} \emph{Let $G$ be a positive discrete unbounded operator on $\H_A$ (Def.\ref{D-H}). A map $\Phi$ in $\Y_G(A,B)$ can be represented as $\Phi=\Psi_1-\Psi_2$, where $\Psi_1$ and $\Psi_2$  are maps in $\F^0_G(A,B)$,  if and  only if there exist linear maps $\Lambda_1$ and $\Lambda_2$ from $\T_G(\H_A)$ into $\T(\H_B)$ such that the map
  $$
  \Theta(\rho)= \left[\begin{array}{cc}
        \Lambda_1(\rho) & \Phi(\rho)\\
        \Phi(\rho) & \Lambda_2(\rho)
        \end{array}\right]
  $$
belongs to the cone $\,\Y_G^{+}(A, B_2)$, where $\H_{B_2}\doteq\H_B\oplus \H_B$.}
\end{property}\medskip

\emph{Proof.} By Theorem \ref{third}B  $\,\Y^+_G(A, B_2)=\F^0_G(A, B_2)$. So, if $\Theta\in\Y^+_G(A, B_2)$ then
there exist a  Hilbert space $\H_E$
and an operator  $V$  in $\B^0_G(\H_A, \H_E\otimes(\H^1_B\oplus\H^2_B))$ such that
$$
 \Theta(\rho)=\Tr_E V\rho V^*, \;\rho\in\T_G(\H_A).
$$
Denote by $P_1$ and $P_2$ the projectors onto the first and the second summands in $\H^1_B\oplus\H^2_B$ correspondingly. Then
the operators $V_1=(P_1\otimes I_E) V$ and $V_2=(P_2\otimes I_E) V$ belong, respectively,  to the spaces $\B^0_G(\H_A, \H^1_B\otimes\H_E)$
and $\B^0_G(\H_A, \H^2_B\otimes\H_E)$. By using these operators one can represent the map $\Theta$ as follows
$$
\Theta(\rho)= \left[\begin{array}{ll}
        \Tr_E V_1\rho V_1^* & \Tr_E V_1\rho V_2^*\\
        \Tr_E V_2\rho V_1^* & \Tr_E V_2\rho V_2^*
        \end{array}\right].
$$
Hence $\Phi(\rho)=\Tr_E V_1\rho V_2^*$. Since the map $\Phi$ is Hermitian-preserving, it follows that $\Phi=\frac{1}{4}(\Psi_+-\Psi_-)$, where
$$
\Psi_+(\rho)=\Tr_E (V_1+V_2)\shs\rho\shs(V_1+V_2)^*\quad  \textrm{and}\quad \Psi_-(\rho)=\Tr_E (V_1-V_2)\shs\rho\shs(V_1-V_2)^*
$$
are maps from the cone $\F^0_G(A,B)$. \smallskip

If $\Phi=\Psi_1-\Psi_2$  then it is easy to see that the map
  $$
\Theta(\rho)=\left[\begin{array}{cc}
        (\Psi_1+\Psi_2)(\rho) & \Phi(\rho)\\
        \Phi(\rho) & (\Psi_1+\Psi_2)(\rho)
        \end{array}\right]
  $$
belongs to the cone $\Y^+_G(A, B_2)$.  $\square$ \smallskip

The above  conjecture is indirectly supported by the following proposition, since
the property stated therein  holds for any map $\Phi$ belonging to the linear span of $\F^0_G(A,B)$
by Proposition \ref{ch-prop} in Section 4.\smallskip

\begin{property}\label{c-nc} \emph{If a map $\,\Phi$ belongs to the completion of the space  $\,\Y(A,B)$ w.r.t.
the ECD norm then $\,\|\Phi\|_{\diamond,E}^G=o\shs(E)\;$ as $\;E\rightarrow+\infty$}.
\end{property}\medskip

\emph{Proof.} There is a sequence $\{\Phi_n\}$ of maps in $\Y(A,B)$ such that $\|\Phi_n-\Phi\|_{\diamond,E_0}^G$ tends to zero as $n\rightarrow+\infty$
any given $E_0>0$.
By using  the triangle inequality and relations (\ref{ECD-n-eq}) we obtain
$$
\left|\frac{\|\Phi_n\|_{\diamond,E}^G}{E}-\frac{\|\Phi\|_{\diamond,E}^G}{E}\right|\leq\frac{\|\Phi_n-\Phi\|_{\diamond,E}^G}{E}
\leq\frac{\|\Phi_n-\Phi\|_{\diamond,E_0}^G}{E_0}
$$
for arbitrary $E>E_0$. Since $\|\Phi_n\|_{\diamond,E}^G=o\shs(E)\,$ as $\,E\rightarrow+\infty\,$ for all $n$ (in  fact, all the functions $E\mapsto\|\Phi_n\|_{\diamond,E}^G$ are bounded), the above inequality implies that  $\|\Phi\|_{\diamond,E}^G=o\shs(E)\,$ as $\,E\rightarrow+\infty$. $\square$\smallskip

\bigskip

I am grateful to A.S.Holevo, G.G.Amosov, A.V.Bulinsky, S.Pirandola, T.Shulman and M.Wilde for discussion and useful remarks.
Special thanks to  Frederik vom Ende for pointing the relation between Theorem 1 in \cite{Ende} and Corollary \ref{main-c} in this paper.

\end{document}